\documentclass{amsart}
\usepackage{amssymb}
\usepackage{graphics}

\theoremstyle{plain}
\newtheorem{lemma}{Lemma}
\newtheorem{proposition}{Proposition}
\newtheorem{theorem}{Theorem}
\theoremstyle{definition}
\newtheorem{definition}{Definition}
\theoremstyle{corollary}
\newtheorem{corollary}{Corollary}
\theoremstyle{definition}
\newtheorem{example}{Example}
\newtheorem{remark}{Remark}
\def \dim {\operatorname{dim}}
\begin{document}
\title[Circle actions on symplectic manifolds]{Semi-free Hamiltonian circle actions on 6 dimensional symplectic manifolds}
\author{Hui Li}
\address{Mathematics Department \\
         University of Illinois \\
         Urbana-Champaign, IL 61801 }
\email{hli@math.uiuc.edu}

\subjclass{Primary : Symplectic geometry; Secondary : Algebraic topology}
\begin{abstract} 
   Assume $(M, \omega)$ is a connected, compact  6 dimensional symplectic manifold equipped with a semi-free Hamiltonian circle
    action, such that the fixed point set consists of isolated points or compact orientable surfaces. We restrict attention to the case $\dim$ $H^2(M)<3$. We give a complete list of the possible  manifolds, determine their equivariant
cohomology ring and equivariant Chern classes. Some of these manifolds are  classified up to diffeomorphism. We also show the existence 
for a few cases. 
\end{abstract}

\maketitle

\section{Introduction} 
    Assume $(M, \omega)$ is a connected, compact 6 dimensional symplectic manifold equipped with a nontrivial semi-free (free outside fixed point sets) Hamiltonian circle
    action. It is well known that the moment map $\phi$ is a real-valued perfect Bott-Morse function. The critical points of $\phi$
    are exactly the fixed point sets $M^{S^1}$ of the circle action, and $M^{S^1}$ is a disjoint union of symplectic submanifolds.
    Therefore, in dimension 6, the fixed point sets could be $0, 2$ or $4$ dimensional. If the fixed points are all isolated, by [TW1], there are precisely
    8 fixed points and the manifold is diffeomorphic to $(\mathbb{C}P^1)^3$. If the fixed point set consists of
 isolated points or  
    orientable closed surfaces, the situation is  much more complicated. In this paper, we consider this question 
    for the case $\dim$ $H^2(M)<3$. We give a complete list of the possible manifolds, determine their equivariant cohomology rings and equivariant Chern classes.
     Some of these manifolds are classified  up to diffeomorphism. One of them arises
as a coadjoint orbit. Some of the others arise as toric varieties.

    Each critical set of the moment map $\phi$  has even index, which is twice of the 
    number of negative weights of the isotropy representation on the normal bundle of the critical set. Since the action is
    semi-free, the weights are $+1$ or $-1$. In dimension 6, there are only index $0, 2, 4$ or $6$ critical sets. We know that
    $\phi$ has a unique local minimum and a unique local maximum ([A]).\\    
    
  \begin{example}\label{ex1}
 Let $S^1$ act on $\mathbb{C}P^3$ as: $\lambda [z_1,z_2,z_3,z_4]=[\lambda z_1,\lambda z_2,z_3,z_4]$.
  The action is semi-free and Hamiltonian with 2 fixed spheres; they are the minimum and the maximum of the moment map.
  \end{example}
  \begin{example}\label{ex2}
   Consider $SO(5)$. Let $T^2$ be a maximal torus. Embed
 $\frak t^*$ in $\frak g^*$ (dual lie algebra of $SO(5)$). Identify $\frak t^*$  with $\mathbb{R}^2$.
The roots of the $SO(5)$ action on $\frak t^*$ are $(\pm 1,0), (0,\pm 1),(\pm 1,\pm 1).$  Recall that the coadjoint orbit through an element $x$ of $\frak t^*$ is a symplectic manifold and the projection to $\frak t^*$ is a  moment map for the $T^2$-action. The fixed points for the $T^2$ action are exactly the Weyl group orbit of $x$ in $\frak t^*$. The isotropy weights at a fixed point $y$ are exactly those roots $\alpha\in \frak t^*$ for which $\langle\alpha,y\rangle <0$. Now consider the coadjoint orbit of $SO(5)$ through the point $(1,0)$. The Weyl group orbit of this point consists of the points $(1,0),(-1,0),(0,1),$ and $(0,-1)$. See the following picture for the moment map image of $T^2$.

    \begin{figure}[h!]
    \scalebox{.80}{\includegraphics{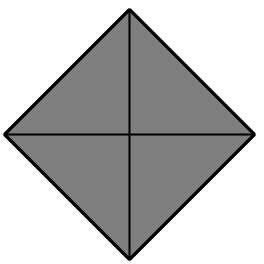}}
    \end{figure}

 The isotropy weights at $(1,0)$ are $(-1,1),(-1,0)$, and $(-1,-1)$. It is easy to see that the $e\times S^1$ action is semi-free and it has 3 fixed point set components: an isolated maximum, an
index 2 sphere and an isolated minimum.\\

  \end{example}

   Let $S^1$ act on a connected, compact symplectic manifold $(M, \omega)$ with moment map $\phi$. Let $a\in$ im$(\phi)$, then $\phi^{-1}(a)/S^1=M_a$ 
   is called the reduced space at $a$.  If the action is semi-free, then for regular value $a$, $M_a$ is a smooth manifold; 
   for singular value $a$, it is not clear whether or not $M_a$ is a smooth manifold. 
   However, in dimension 6, the fixed point sets are the minimum, the maximum or of index 2 or co-index 2.  By [GS], all the reduced spaces 
   (including those on critical levels) are smooth manifolds. Moreover, when the index 2 and co-index 2
fixed point set consists of surfaces, the reduced spaces are all diffeomorphic, and the index 2 (also of co-index 2)  surface is symplectically
 embedded in this reduced space.\\

    Now assume the fixed point set only consists of surfaces. Then the reduced spaces are all diffeomorphic to
 an $S^2$ bundle over a Riemann surface. (It is enough to consider the reduced space at a regular value right
 above the minimum Riemann surface). Moreover, the natural diffeomorphism between the reduced spaces is induced by the Morse flow on $M$.
 More explicitly, if there is no critical set between two level sets, the Morse flow gives a one to one correspondence between the level sets.
  Assume there is  an index 2 (also of co-index 2) fixed surface $F$ with $\phi(F)=c$. Assume $F$ is the only fixed point set in $\phi^{-1}(c-\epsilon, c+\epsilon)$.
  Then each $S^1$ orbit in this neighborhood which flows into $F$ corresponds to a unique point of $F$ and corresponds
 to a unique $S^1$ orbit which flows out of $F$. (See [Mc] for detail). So in both cases, there is a well defined diffeomorphism between the reduced
 spaces. This diffeomorphism
  takes a homology class of the reduced space at a regular level right above the minimum to a homology class of
 the reduced space at a regular level right below the maximum. We give the following definition:
 \begin{definition}
   Let $(M,\omega)$ be a compact connected symplectic manifold equipped with a semi-free Hamiltonian circle action.
  Assume the fixed point set consists only of surfaces. If the fibers of the reduced spaces at a regular level right
 above the minimum and right below the maximum can be represented by the same cohomology class of the reduced space,
 we say there is $\bold{no\, twist}$; otherwise, we say there is a $\bold{twist}$.
 \end{definition}

  $\mathbb{C}P^3$ in Example~\ref{ex1} has a twist, see Section 6 for a picture and an explanation.
 \begin{example}\label{ex3}
  Let $S^1$ act on ${(\mathbb{C}P^1)}^3$ by rotating the first 2 spheres with speed
 $1$ and fixing the third sphere; the fixed point set consists of 4 spheres, and there is no twist.       
 \end{example}   
 
   \begin{definition}
  For each  manifold, we  define the $\bold{fixed\, point\, data}$ to be the diffeomorphism type of the fixed point sets, their indices,
 and $b_{min}$, where $b_{min}$ denotes the 1st Chern number of the normal bundle of the minimum when the minimum is a surface.
\end{definition}

\begin{remark}\label{rem0}
 The fixed point data is not enough to describe the manifold type. In Section 6, we give two toric varieties  which have the same fixed point data but one is twisted and the other one is not twisted. 
\end{remark}

    \begin{theorem}
     Let $(M,\omega)$ be a  connected compact 6 dimensional symplectic manifold equipped with a semi-free Hamiltonian circle action such that the fixed point set consists of isolated  points or  surfaces.  Assume $\dim$ $H^2(M)<3$. Then the possible fixed point data and the twist type are:

    \begin{figure}[h!]
    \scalebox{.60}{\includegraphics{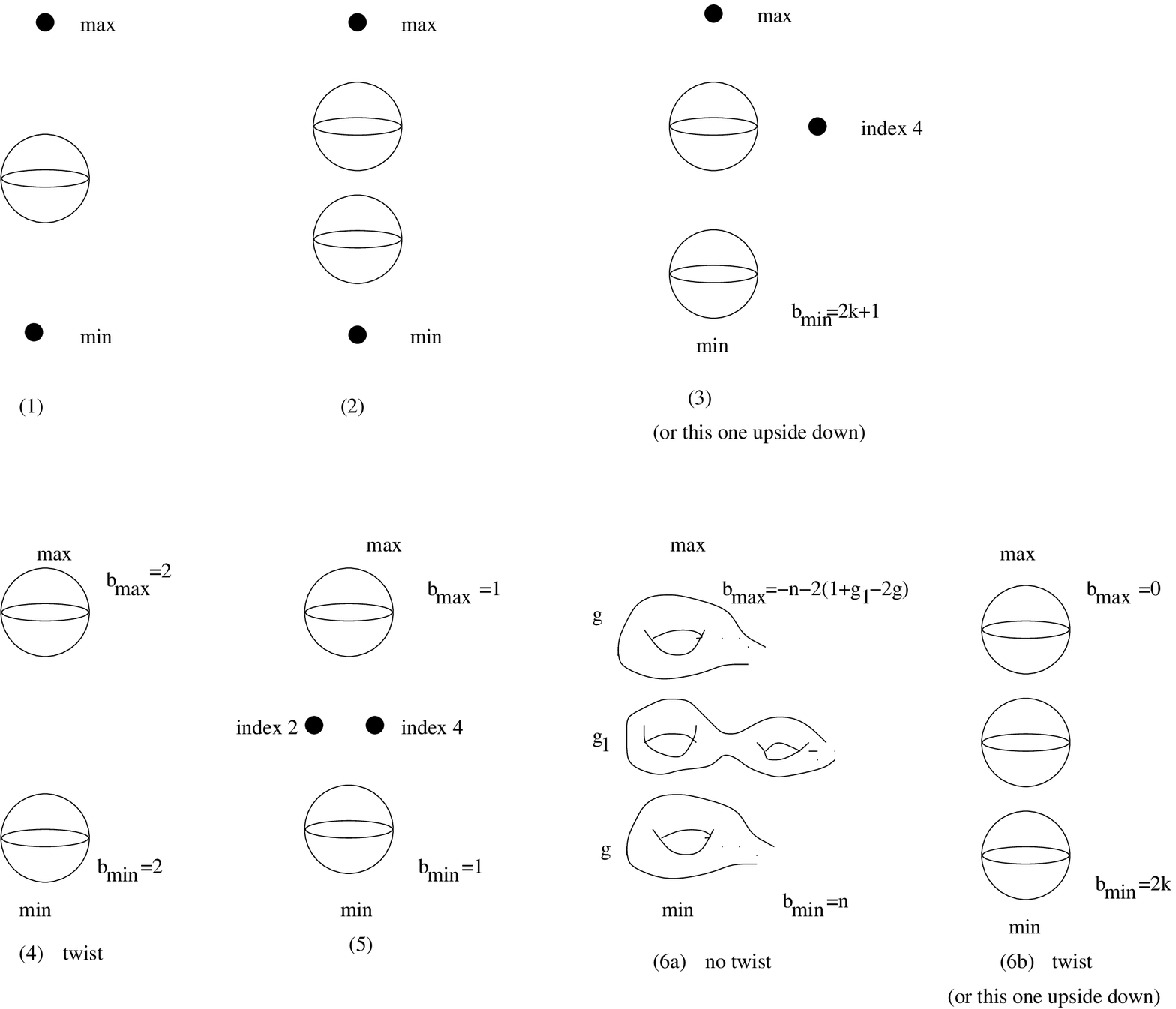}}
    \end{figure}

 In the picture, $g$, $g_1$ are the genus of the surfaces, and $b_{max}$ is the 1st Chern number of the normal bundle of the  maximum.
    \end{theorem} 

\begin{remark}\label{rem1}
 The picture gives  us an order of how the moment map crosses the fixed point sets. For instance, manifold of type (1) means that the manifold has an isolated minimum,
an index 2 sphere, and an isolated maximum.

    For manifold of type (2), the two index 2 spheres can not be on the same level of the moment map.

    For manifolds of type (3) and (5), the moment map can cross the two non-extremal fixed point sets  in any order or the 2 fixed point sets are on the same level of the moment map.
\end{remark}

   Notice that Example 1 is of type (4)  and  Example 2 is of type (1), so  this proves the existence of manifold of type (1) and type (4).

  \begin{proposition}
  Let $b_+$ and $b_-$ be the Euler numbers of the positive normal bundle and the negative  normal bundle of an index 2 surface.

  Then  for each  manifold which has an index 2 surface in Theorem 1, $b_+$ and $b_-$ of each index 2 fixed surface are determined by the fixed point data and the twist type.
\end{proposition}
  \begin{theorem}
     Let $(M,\omega)$ be a  connected compact 6 dimensional symplectic manifold equipped with a semi-free Hamiltonian circle action such that the fixed point set consists of isolated  points or  surfaces.  Assume $\dim$ $H^2(M)<3$. Then the equivariant cohomology ring of $M$ and its restriction to the fixed point
 set are uniquely determined by the fixed point data and the twist type.
 \end{theorem}

 \begin{corollary}
  For each manifold in Theorem 1, the equivariant Chern classes and their restrictions to the fixed point set are uniquely determined by the fixed point data and the twist type.
 \end{corollary}  
  \begin{proof}
 The restrictions of the equivariant Chern classes to the fixed point sets are determined by the fixed
 point sets, their indices and the Euler numbers of the normal bundles of the fixed point sets. Therefore Proposition 1 and Theorem 2 
 imply the corollary.  
  \end{proof}

 \begin{corollary}\label{cor2}
  For each manifold in Theorem 1, the ordinary cohomology ring and the ordinary Chern classes are uniquely determined 
 by the fixed point data and the twist type.
 \end{corollary}

  \begin{proof}
  By Theorem 2 and Corollary 1, we only need to show that  the map $H^*_{S^1}(M, \mathbb{Z})\rightarrow H^*(M, \mathbb{Z})$ is
 surjective\footnote{This argument is suggested by Franz Matthias.}. For each fixed point set component $F$ such that $\phi(F)=c$, assume $F$
 is the only critical set  in  $\phi^{-1}(c-\epsilon, c+\epsilon)$. Let $M_-=\{x\in M |\, \phi(x)<c-\epsilon\}$ and $M_+=\{x\in M |\, \phi(x)<c+\epsilon\}$.
  By induction, assume $H^*_{S^1}(M_-,\mathbb{Z})\rightarrow H^*(M_-, \mathbb{Z})$ is surjective. We have the following commutative diagram:
  \[
   \begin{array}{ccccccccc}
   0 &\rightarrow & H^*_{S^1}(M_+, M_-,\mathbb{Z})& \rightarrow &H^*_{S^1}(M_+,\mathbb{Z}) &\overset{f}\rightarrow &H^*_{S^1}(M_-,\mathbb{Z}) &\rightarrow &0\\
     &           &\downarrow h_1                      & &\downarrow h_2 &              & \downarrow h_3       &   &     \\
   0 &\rightarrow &H^*(M_+, M_-, \mathbb{Z})&\rightarrow       &H^*(M_+,\mathbb{Z})& \overset{g}\rightarrow   &     H^*(M_-,\mathbb{Z})& \rightarrow &0\\
  \end{array}
  \]
  
   We have the commutative diagram for the long exact sequences in equivariant cohomology and ordinary cohomology. 
   The long exact sequence for equivariant cohomology  in $\mathbb{Z}$ coefficients splits is by [TW2] (since the fixed point set has no torsion 
   cohomology). The surjectivity of $H^*_{S^1}(M_+, \mathbb{Z})\rightarrow H^*_{S^1}(M_-, \mathbb{Z})$ and the surjectivity of  
  $H^*_{S^1}(M_-,\mathbb{Z})\rightarrow H^*(M_-, \mathbb{Z})$ give the surjectivity of $H^*(M_+, \mathbb{Z})\rightarrow H^*(M_-, \mathbb{Z})$ (for any
  $*$). Therefore we have the short exact sequence for ordinary cohomology.  $h_1$ is clearly surjective after
  we identify the relative equivariant  cohomology (or cohomology) with the equivariant cohomology (or cohomology) of $F$ by using Thom isomorphism.
      By diagram chasing, we can prove $h_2$ is surjective.   
        
\end{proof}

   Wall classified  certain  6-dimensional manifolds. For symplectic manifolds, his theorem implies the following 

  \begin{theorem} ([W])
  If two simply-connected 6 dimensional symplectic manifolds have  no torsion homology, their second
        Stiefel-Whitney classes $w_2$ vanish, and they have the same cohomology ring and Chern classes, then they are diffeomorphic.
  \end{theorem}

 We can use this to  prove

 \begin{proposition}
  All manifolds of type (4) are diffeomorphic to $\mathbb{C}P^3$. Hence it is diffeomorphic to a toric variety.

  All manifolds of type (6a) with $g=0$, the same $g_1$ and the same even $b_{min}$ are diffeomorphic.

  All manifolds of type (6b) with the same $b_{min}$   are diffeomorphic. 
 \end{proposition}

 We will also prove 
\begin{proposition}
 If  $(M,\omega)$ is of type (6b) and $b_{min}=2$, then $(M,\omega)$ is diffeomorphic to a toric variety which is 
 $\mathbb{C}P^3$ blown up a point.
  
 There exists a toric variety  when  $(M,\omega)$ is of type (3) and $b_{min}=1$ or $b_{min}=3$.
 In particular, when $b_{min}=1$, the toric variety is $\mathbb{C}P^2\times\mathbb{C}P^1$. 
\end{proposition}

   \subsubsection*{Acknowledgement}
 I would like to give  special thanks to Susan Tolman for many helpful discussions, especially for explaining to me the integration formula in
 equivariant cohomology and the Delzant polytope.
She also gave me many good suggestions about how to write this paper.

   I would like to thank Miguel Abreu for suggesting to construct the polytopes. The polytopes helped my understanding of the problem and
   provided good examples.

 I also would like
 to thank Jeremy Wang  for patiently teaching me how to draw pictures.

   Finally, I would like to thank the referee for many constructive comments.

\section{Proof of Theorem 1}

\begin{lemma} \label{lem1}  When the maximum and the minimum are isolated or they are surfaces, if there are  isolated fixed points with index 2 or index 4, then they must
 occur in pairs.

 When the minimum is a surface, and the maximum is a point, then the minimum is a sphere. If  index 2 or index 4 isolated points occur, the number of index 4 points
 is $1$ more than the number of index 2 points.
\end{lemma}
             
\begin{proof}

   Use the fact that $\phi$ is a perfect  Morse-Bott function  and Poincare duality. When the minimum and the maximum are isolated, they don't contribute to  $H^2(M)$ or $H^4(M)$. Each
     index 2 surface increases both $\dim$ $H^2(M)$ and $\dim$ $H^4(M)$ by $1$. Each index 2 isolated fixed point increases 
     $\dim$ $H^2(M)$ by $1$, $\dim$ $H^4(M)$ should also be increased by $1$ which needs the index 4 isolated point.
     
       When the minimum and the maximum are surfaces, the minimum increases $\dim$ $H^2(M)$ by $1$, and  the maximum increases $\dim$ $H^4(M)$ by $1$. Reasoning  as
      above, we can prove the corresponding statement.

   If the minimum is a surface, and the maximum is isolated, then the minimum has to be a sphere. The reason is: $\phi$ is a perfect Morse-Bott function and each critical set has even index, so none of the critical sets contributes to $H^5(M)$,
 hence $H^5(M)=0$; by duality, $H^1(M)=0$. Since to pass these
 critical sets won't increase or kill $H^1(M)$, therefore the minimum surface has genus 0. Argue similarly for the rest of the statement.
       
  \end{proof}
 
  We will approach the proof of Theorem 1 separately, according to the dimensions of the minimum and the maximum.
 \subsection{The minimum and the maximum are isolated}

 If $\dim$ $H^2(M)<3$, by Lemma~\ref{lem1}, we have the following 3 possibilities  
       ($N_F$ denotes the number of fixed point set components):

       (i) $N_F=3$: the minimum,  one index 2 surface,  and the maximum.
       
       (ii) $N_F=4$: the minimum,  two index 2 surfaces, and the maximum.
       
       (iii) $N_F=5$: the minimum, one index 2 surface, one index 2 point,  one  index 4 point, and the maximum.
      
       ( There is at least 1 surface, see introduction for reason.)
       
  In fact the following is true:
  \begin{lemma} \label {lem2}    
       Let $(M, \omega)$ be a connected, compact 6 dimensional symplectic manifold equipped with a semi-free Hamiltonian circle action
       such that the fixed point set consists of isolated points or surfaces. If the maximum and the minimum are isolated, and the other fixed point
        sets are surfaces, then (i) and (ii) are the only possibilities
        with the surfaces being spheres. 
   \end{lemma}

         Before proving the lemma, let's do some preparation.

      If $a$ is a regular value of the moment map $\phi$, then $S^1$ acts freely on $\phi^{-1}(a)$, so we have a principal circle bundle

\begin{equation}\label{eq1}    
\begin{array}{ccl}
        S^1 & \hookrightarrow & \phi^{-1}(a)=P_a    \\                           
            &  & \downarrow    \\          
                           & &   \phi^{-1}(a)/S^1=M_a, \\
 \end{array}                          
  \end{equation}                             
where $M_a$ is the reduced space at $a$.
           
             Let $c$ be a critical value of the moment map $\phi$, and  $S\subset \phi^{-1}(c)$ be an index 2 fixed surface.  Assume for small $\epsilon$, $S$ is the only critical set in $\phi^{-1}([c-\epsilon, c+\epsilon])$. Let $P_{c-\epsilon}$ be the circle bundle $\phi^{-1}(c-\epsilon)$ over $M_{c-\epsilon}$, $e(P_{c-\epsilon})$ be its Euler class.
              $P_{c+\epsilon}$ and $e(P_{c+\epsilon})$ have similar meanings. By [GS] or [Mc], 
             $M_{c+\epsilon}$ is diffeomorphic to $M_{c-\epsilon}$ (in dimension 6), let's call them $M_{red}$ (reduced space). Moreover, there is a symplectic embedding
             $i: S\rightarrow M_{red}$ such that  if $i(S)=Z$, and $\eta$ is the dual class of $Z$ in 
             $M_{red}$, then
 \begin{equation} \label{eq2}           
             e(P_{c+\epsilon})=e(P_{c-\epsilon})+\eta.                  
 \end{equation}

             Choose an almost complex structure of $M_{red}$ such that $Z$ is J-holomorphically embedded in $M_{red}$, then we have the adjunction formula 
            (see for example [McS])
      \begin{equation} \label{eq3}           
                Z\cdot Z-<c_1(TM_{red}), Z>+2=2g,             
   \end{equation}             
              where $Z\cdot Z$ is the self-intersection of $Z$ in $M_{red}$, $c_1(TM_{red})$ is the 1st Chern class of the tangent bundle of $M_{red}$,
           $<,>$ is the natural pairing between homology and cohomology, and $g$ is the genus of the surface $S$.\\

     In this article, when we say ``above" or ``below" some critical level, we mean at a regular level in a small neighborhood of the critical
      level, such that there are no other critical level in this neighborhood.

     Proof of Lemma~\ref{lem2}:
             \begin{proof}                          
                Since the minimum (or maximum)is isolated, $M_{red}$ at a regular value above the minimum (or below the maximum) is diffeomorphic to $\mathbb{C}P^2$, 
                and the circle bundle (\ref{eq1}) over $M_{red}$ is the Hopf fibration
     \[                    
  \begin{array}{ccl} 
               
                   S^1 & \hookrightarrow & S^5\\
                         &                 &  \downarrow \\
                         &                 &  \mathbb{C}P^2.
                       
  \end{array}
       \]                
                                          
               If $u\in H^2(\mathbb{C}P^2)$ is a generator, and$-u$ is the Euler class of this circle bundle at regular values above the minimum,
               then $u$ is the Euler class of it at regular values  below the maximum. If $\eta=au$ is the dual class of a surface in $\mathbb{C}P^2$,
               then by (\ref{eq3})
           $$\int_{\mathbb{C}P^2}(au)^2-\int_{\mathbb{C}P^2}3u \cdot au +2=2g,$$ 
    i.e.,
    \begin{equation} \label{eq4}
            a^2-3a+2=2g.                                
    \end{equation}

    If $a_1u, a_2u,...$ are the dual classes of the fixed surfaces in $\mathbb{C}P^2$, then by (\ref{eq2}), we have
    \begin{equation} 
               -u+a_1u+a_2u+...=u.
    \end{equation}
     The surfaces are symplectically embedded in $\mathbb{C}P^2$ which implies $a_i>0$ for $i=1, 2,...$,
     so either $a=2$, i.e., there is only one fixed surface with dual class  $\eta=2u$ in $\mathbb{C}P^2$; or $a_1=1, a_2=1$, i.e., 
     there are two fixed surfaces with dual classes
     $\eta_1=\eta_2=u$ in $\mathbb{C}P^2$.  In either of these two cases, (\ref{eq4}) gives 
     $g=0$.

  Notice that for the case with 2 spheres, the 2 spheres can't be on the same level of the moment map. If they were, they would have intersection 
      $\int_{\mathbb{C}P^2}u\cdot u=1$ which contradicts  the fact that they are disjoint.

    \end{proof}

Notice from the proof the following. Assume the reduced spaces are all diffeomorphic to $\mathbb{C}P^2$. If there
is any fixed point, then the action has to be Hamiltonian.\\

    Regarding (iii), let's prove
    \begin{lemma}\label{lem3}
     There are no manifolds in category (iii).
     \end{lemma}

    By [GS], when the moment map crosses an index 2 point (index 4 point), the reduced space changes by a blow up at a point (blow down of an exceptional
    divisor). Now let $M_-$ be the reduced space before surgery (blow up or blow down), $M_+$ be the reduced space after surgery, and let
    $e(P_-)$ be the Euler class of the principal circle bundle over $M_-$, $e(P_+)$ be the Euler class of the principal circle bundle over $M_+$. 
    Let $\eta_{ex,div}$ be the dual class of the exceptional sphere in $M_+$ after blow up or in $M_-$ before blow down. Let's always call the blow
     down map by $\beta$. By [GS], we have the following:
     
      Corresponding to blow up,
    \begin{equation} \label{eq6}
       e(P_+)=\beta^*e(P_-)+\eta_{ex,div}.
    \end{equation}

      Corresponding to blow down,
    \begin{equation} \label{eq7}
       \beta^*e(P_+)=e(P_-)+\eta_{ex,div}.
    \end{equation}
    \begin{remark}
      These are the formulas for $\phi$ going from minimal value to maximal value. When we look at the manifold upside down, 
      the Euler class of the principal circle bundle changes sign. So the above two formulas don't contradict  each other.
     \end{remark}

    Proof of Lemma~\ref{lem3}:
    \begin{proof}
   Assume the index 2 surface has genus $g$.

     Since the reduced space at a regular value above the minimum is $\mathbb{C}P^2$, there is no exceptional sphere to blow down,
     so the moment map must cross the index 2 point before it crosses the index 4 point. The moment map crosses the critical sets in 2 possible
     orders: (a) the minimum, the index 2 surface, the index 2 point, the index 4 point, the maximum. (b) the minimum, the index 2 point, the index 2 surface, the index 4 point, the maximum.
     (The index 2 surface next to the maximum case is case  (a) with the action reversed.)
     
     Let $\tilde{\mathbb{C}P^2}$ be  $\mathbb{C}P^2$ blown up at a point.  We look at $\tilde{\mathbb{C}P^2}$ as the nontrivial $S^2$ bundle  over $S^2$.
      Let $x, y\in H^2(\tilde{\mathbb{C}P^2)}$ be a basis which
     are the dual class of the fiber and the dual class of the exceptional divisor.

     For (a), the reduced spaces change in the following way:
     \begin{figure}[h!]
    \scalebox{.57}{\includegraphics{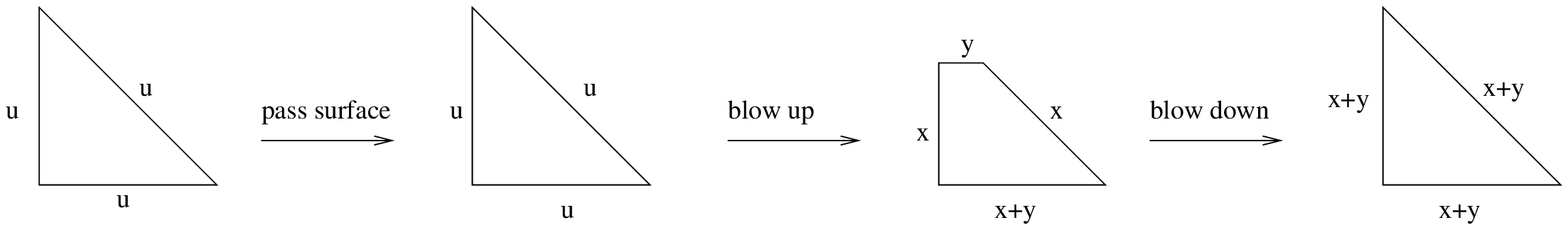}}       
    \caption{Reduced Spaces}\label{labelname:xxx}
    \end{figure}

     Let $au$ be the dual class of the index 2 surface in $\mathbb{C}P^2$, by (\ref{eq2}), (\ref{eq6}) and (\ref{eq7}), we have
     
      $$ -u+au+y+y=x+y,$$  where $\beta^*u=x+y$.

      We see that there is no $a$ satisfying this equation. So (a) is not possible.

      For (b), the reduced spaces change in the following way:
       \begin{figure}[h!]
       \scalebox{0.57}{\includegraphics{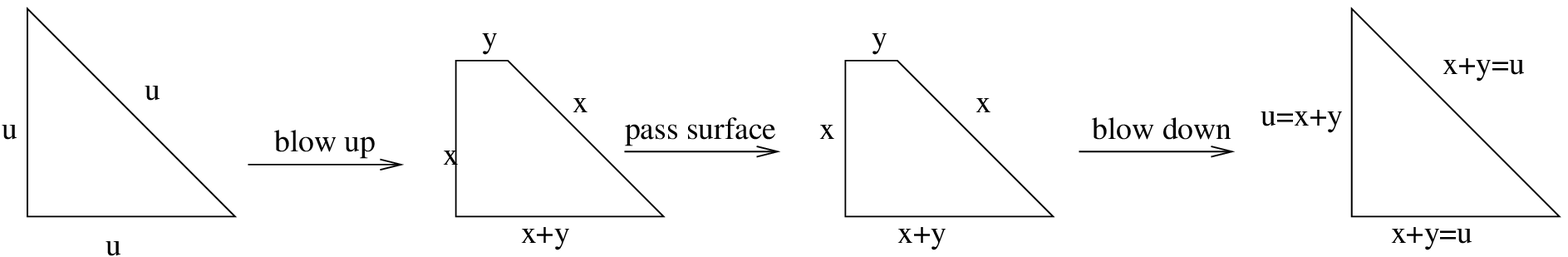}}
       \caption{Reduced Spaces}\label{labelname:2}
       \end{figure}
   
     In this case, the index 2 surface is embedded in $\tilde{\mathbb{C}P^2}$. Let $cx+dy$ be the dual class of the index 2 surface in 
     $\tilde{\mathbb{C}P^2}$, then by (\ref{eq2}), (\ref{eq6}) and (\ref{eq7}),
      we have
     $$-u+y+(cx+dy)+y=x+y,$$  where $\beta^*u=x+y$.
    we get
      $$c=2, \quad d=0.$$
         
       The adjunction formula (\ref{eq3}) is
       
        $$\int_{\tilde{\mathbb{C}P^2}}(cx+dy)^2-\int_{\tilde{\mathbb{C}P^2}}(3x+2y)(cx+dy)+2=2g,$$
        i.e.,
       \begin{equation}\label{eq8}
        2cd-d^2-d-2c+2=2g.
        \end{equation}
       For $c=2, d=0$, there is no $g\geqslant 0$  satisfying this formula. So (b) is not possible either.
       
       From the above, we know that this manifold doesn't exist.
     \end{proof}

     So in this section, we proved there are only 2 possible manifolds (i) and (ii) with the surfaces being spheres. Let us
     give them new names manifold of type (1) and manifold of type (2). They are type (1) and type (2) in Theorem 1.

     \subsection{The minimum is a surface and the maximum is a point}

      By Lemma~\ref{lem1}, the surface is a sphere.

    \begin{lemma} \label{lem4}
    Assume the maximum is isolated, the minimum is a sphere, and there are no  index 2 surfaces. Let $N_2$ and $N_4$  denote the number of index 2 and index 4 
    isolated fixed points respectively, then $N_4=3, N_2=2$.
    \end{lemma}
    
     Thus, $\dim$ $H^2(M)\geq 3$, and we don't consider this case.
     
    \begin{proof}
     We will use integration formula in equivariant cohomology. (See [Au] or [AB] for this formula.)
     Let $\lambda\in H^*_{S^1}(pt)$ be a generator, and $u\in H^2(\hbox{minimum})$ be a generator. Let $b_{min}$ be the 1st Chern number of the normal
     bundle of the minimum.
    
    Integrate $1$ on $M$:
      \[
         \int_{\hbox{min}}\frac{1}{\lambda^2+b_{min}\lambda u}+\frac{1}{(-\lambda)^3}|_{\hbox{max}}+\sum_{\hbox{index 2 points}}\frac{1}{-\lambda\cdot \lambda\cdot \lambda}+\sum_{\hbox{index 4 points}}\frac{1}{-\lambda\cdot 
          -\lambda\cdot \lambda}=0.
      \]
   We get  \qquad         $-b_{min}-1-N_2+N_4=0.$\\
        
     Consider the equivariant 1st Chern class  $c_1(TM)$:\\
        $c_1(TM)|_{\hbox{min}}=2\lambda+2u+b_{min}u$.\\
       $c_1(TM)|_{\mbox{index 2 point}}=-\lambda+\lambda+\lambda=\lambda$.\\
       $c_1(TM)|_{\mbox{index 4 point}}=-\lambda-\lambda+\lambda=-\lambda$.\\
       $c_1(TM)|_{\hbox{max}}=-\lambda-\lambda-\lambda=-3\lambda$.\\

     Integrate $c_1(TM)$ on $M$:

       \[
          \int_{\hbox{min}}\frac{2\lambda+2u+b_{min}u}{\lambda^2+b_{min}\lambda u}+\frac{-3\lambda}{(-\lambda)^3}+\sum_{\hbox{index 2 points}}\frac{\lambda}{-\lambda^3}+\sum_{\hbox{index 4 points}}\frac{-\lambda}{\lambda^3}=0.
       \]   
     We obtain   \qquad $-b_{min}+5-N_2-N_4=0.$ Combine the above two equations, we get $N_4=3$. By Lemma~\ref{lem1}, $N_2=2$.
     \end{proof}

     If we assume $\dim$ $H^2(M)<3$, then by Lemma~\ref{lem4} and Lemma~\ref{lem1}, we  have the following:\\
    (iv) isolated maximum, sphere minimum, an index 4 point, and an index 2 surface.

     \begin{lemma}\label{lem5}
      The index 2 surface in (iv) is a sphere.
     \end{lemma}
     
      These are manifolds of type (3) in Theorem 1.

     Before proving Lemma~\ref{lem5}, let's look at the reduced space at a regular value of $\phi$ above the minimum sphere. It is an $S^2$ bundle over $S^2$, so it is either diffeomorphic
     to $S^2\times S^2$ or diffeomorphic to $\tilde{\mathbb{C}P^2}$ ($\mathbb{C}P^2$ blown up at a point).\\

     Let's analyze this bundle more generally. Assume the minimum is a Riemann surface $\Sigma_g$. The reduced space at a regular value above the minimum is
      an $S^2$ bundle over $\Sigma_g$.

     For $S^2\times\Sigma_g$, let $x, y\in H^2(S^2\times\Sigma_g)$ be a basis which are the dual classes of the fiber $S^2$ 
     and the base $\Sigma_g$ respectively. Then $\int_{base}x=\int_{S^2\times\Sigma_g}xy=1, \quad  \int_{fiber}y=\int_{S^2\times\Sigma_g}yx=1, \quad 
      \int_{fiber}x=\int_{S^2\times\Sigma_g}x^2=0, \quad \int_{base}y=\int_{S^2\times\Sigma_g}y^2=0$.

       Similarly, for the non-trivial bundle $E_{\Sigma_g}$, let $x, y\in H^2(E_{\Sigma_g})$ be a basis which are the dual classes of the fiber and the 
      section $\Sigma_-$ which has self-intersection $-1$ respectively. Then $\int_{\Sigma_-}x=\int_{E_{\Sigma_g}}xy=1, \quad  \int_{fiber}y=\int_{E_{\Sigma_g}}xy=1, \quad 
       \int_{fiber}x=\int_{E_{\Sigma_g}}x^2=0, \quad   \int_{\Sigma_-}y=\int_{E_{\Sigma_g}}y^2=-1$.\\

       Now symplectically cut the manifold ([L]) at a regular value $a$ of $\phi$ above the minimum surface (in a small neighborhood, such that 
      there are no other critical sets). We get a space whose interior is diffeomorphic to $\phi^{-1}(\hbox{minimal}, a)$, whose boundary is $M_a$.
       This space has a Hamiltonian circle action such that 
       it has the same surface as the minimum and it has $M_a$ as the maximum.
      The surface minimum has the same normal bundle as it does in $M$, and the maximum  has  normal bundle $-P_a$.\\
      \indent  Let $b_{min}$ be the 1st Chern number of the normal bundle of the minimum surface, and let $e(P_a)$ be the ordinary Euler class of 
      $P_a$. Integrate $1$ on the cut  space:
      
    $$\int_{M_a}\frac{1}{-\lambda-e(P_a)}+\int_{\hbox{min}}\frac{1}{\lambda^2+b_{min}\lambda u}=0.$$
    We have
                  
      \begin{equation} \label{eq9}     
         \int_{M_a}e^2(P_a)+b_{min}=0.
      \end{equation}

      Let's assume $e(P_a)=px+qy$.
      Since the restriction of the circle bundle (\ref{eq1}) to the fiber is the Hopf fibration
               
     \[         
     \begin{array}{ccl}          
             
          S^1  & \hookrightarrow &  S^3 \\
           &  &  \downarrow\\
           &   &     S^2,\\

     \end{array}
     \]
   so $\int_{fiber}e(P_a)=\int_{M_a}(px+qy)x=q=-1.$

     \indent By (\ref{eq9}), if $M_a\simeq S^2\times\Sigma_g$, then $p=\frac{b_{min}}{2}$, so $b_{min}$ is even. If $M_a\simeq E_{\Sigma_g}$,
      then $p=\frac{b_{min}-1}{2}$, so $b_{min}$ is odd. So we have proved the following:
      \begin{lemma} \label{lem6}
       If the minimum is a surface $\Sigma_g$ with 1st Chern number $b_{min}$ of it's normal bundle, then $M_a$ (a above the minimum) is diffeomorphic to $S^2\times
       \Sigma_g$ if and only if $b_{min}=2k$ is even, and it is diffeomorphic to $E_{\Sigma_g}$ if and only if $b_{min}=2k+1$ is odd.
       In either case, $e(P_a)=kx-y$.
       \end{lemma}
       
     We can prove similarly the following:
     \begin{lemma}\label{lem6'}
     If the maximum is a surface $\Sigma_g$ with 1st Chern number $b_{max}$ of it's normal bundle, then $M_a$ (a below the maximum) is diffeomorphic to 
     $S^2\times\Sigma_g$ if and only if $b_{max}=2k'$ is even, and it is diffeomorphic to $E_{\Sigma_g}$ if and only if 
     $b_{max}=2k'+1$ is odd. In either case, $e(P_a)=-k'x+y$.\\ 
     \end{lemma} 
       
     Proof of Lemma~\ref{lem5}:
     \begin{proof}
     
      In (iv), since $M_{red}$ can only change by a diffeomorphism and a blow down, and $M_{red}$ at a regular value below the maximum is diffeomorphic
       to $\mathbb{C}P^2$, so $M_{red}$ at a regular value above the minimum has to be diffeomorphic to $\tilde{\mathbb{C}P^2}$. Hence we can assume  $b_{min}=2k+1$.
      Let $g$ be the genus of the index 2 surface.

      Assume first that $\phi$  crosses the index 4 point before it crosses the index 2 surface $S$, i.e., $M_{red}$ changes 
      from $\tilde{\mathbb{C}P^2}$ to $\mathbb{C}P^2$ and remains to be $\mathbb{C}P^2$ after
      $\phi$ crosses the index 2 surface.

    Let $x, y$ be the basis of $H^2(\tilde{\mathbb{C}P^2})$ as taken above. After blow down, $x+y\in H^2(\mathbb{C}P^2)$ is the generator. 
     If $\eta=a(x+y)$ is the dual class of the index 2 fixed surface in $\mathbb{C}P^2$,
      then by (\ref{eq2}) and (\ref{eq7}), we have 
           $$(kx-y)+y+a(x+y)=x+y.$$
    So
      
            $$a=1,   \qquad    k=0.$$
      By(\ref{eq4}), if $a=1$, then $g=0$, i.e., the index 2 surface is a sphere. And we  have $b_{min}=1$ in this case.

      Next assume $\phi$ crosses the index 2 surface first, and then the index 4 point. This time the index 2 genus $g$ surface
      is symplectically embedded in $\tilde{\mathbb{C}P^2}$. Let $\eta=cx+dy$ be the dual class of the embedded image of it in $\tilde{\mathbb{C}P^2}$,
       then
       $$(kx-y)+(cx+dy)+y=x+y.$$
      This gives
         $$d=1, \qquad  k+c=1.$$
      If $d=1$, then (\ref{eq8}) gives  $g=0$, i.e., the index 2 surface is a sphere. And the dual class of the index 2 sphere in 
      $\tilde{\mathbb{C}P^2}$ is $\eta=(1-k)x+y$.

       So we get manifold of type (3) with $b_{min}=2k+1$. 
     \end{proof}
     
     \begin{remark}
      Note that if we have manifold(s) of type (3), we also have manifolds with the reverse $S^1$ action.
      \end{remark}
    
    \subsection{The minimum and the maximum are surfaces}

     The reduced space only changes by blowing up a point or blowing down an exceptional sphere, or by a diffeomorphism.
     Neither of these operations will change $\pi_1$ or $H^1$ of the reduced space, so the minimum and the maximum must have the same genus
     in order to keep $\pi_1$ or $H^1$ of the reduced space at regular values above the minimum and below the maximum the same.

     For $\dim$ $H^2(M)<3$, there are 3  possibilities:\\
     (v)The minimum and the maximum with the same genus $g$ and no other fixed point sets.\\
     (vi)The minimum and the maximum with the same genus $g$, 1 index 2 isolated point and 1 index 4 isolated point.\\
     (vii)The minimum and the maximum with the same genus $g$ and an index 2 surface with genus $g_1$.    
     
      About (v), we can prove the following lemma:
     \begin{lemma}\label{lem7}
     If the minimum surface and the maximum surface are the only fixed point sets, then they must be spheres with $b_{min}=b_{max}=2$, and there is a twist
     in this case.
     \end{lemma}
     \begin{proof}
     This follows immediately from Lemma~\ref{lem6} and Lemma~\ref{lem6'}: if there were no twist, then $kx-y=-k'x+y$ which can not be true.
     Only when there is a twist, i.e., $kx-y=-k'y+x$, (v) is possible. Then $k=k'=1$.
     \end{proof}
     This is manifold of type (4) in Theorem 1.

     \begin{lemma}\label{lem8}
      In (vi), $g=0$ and $b_{min}=b_{max}=1$.
     \end{lemma}
     
     This is manifold of type (5) in Theorem 1.
     
     \begin{proof}
     
     Now let's compute the 1st Chern numbers of the normal bundles of the minimum and the maximum:
              
              Integrate $1$ on $M$:

              $$\int_{\hbox{min}}\frac{1}{\lambda^2+b_{min}\lambda u}+\int_{\hbox{max}}\frac{1}{\lambda^2-b_{max}\lambda u}+\frac{1}{-\lambda^3}+\frac{1}{\lambda^3}=0.$$
              We get
              
              $$b_{max}=b_{min}.$$
              
              Integrate the equivariant 1st Chern class $c_1(TM)$ on $M$:
              
              $$\int_{\hbox{min}}\frac{2\lambda+(2-2g)u+b_{min}u}{\lambda^2+b_{min}\lambda u}+\int_{\hbox{max}}\frac{-2\lambda+(2-2g)u+b_{max}u}{\lambda^2-b_{max}\lambda u}+\frac{\lambda}{-\lambda^3}$$
               
              $$+\frac{-\lambda}{\lambda^3}=0.$$
              We get
                   $$b_{min}+b_{max}=2-4g.$$
              
              Combine the above 2 equations, we get
             $$b_{max}=b_{min}=1-2g.$$
          So $M_{red}$ at regular values above the minimum and below the maximum are the nontrivial $S^2$ bundle over the surface.

          If $g=0$, there are 2 possible orders for the moment map to cross the critical sets: (a) the index 4 point first, (b) the index 2 
          point first.

          For (a), the reduced spaces change in the following way:
          \begin{figure}[h!]
       \scalebox{0.60}{\includegraphics{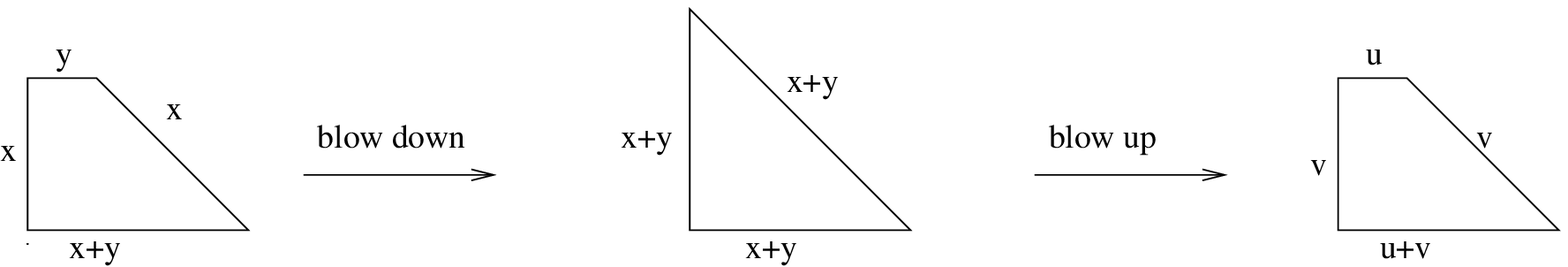}}
       \caption{Reduced Spaces}\label{labelname:3}
       \end{figure}\\
       
      By Lemma~\ref{lem6}, Lemma~\ref{lem6'}, (\ref{eq7}) and (\ref{eq6}),  the Euler classes of the circle bundles change:
           
           $$-y+y+u=u.$$
           
           For (b), the reduced spaces change in the following way:
            \begin{figure}[h!]
      \scalebox{0.60}{\includegraphics{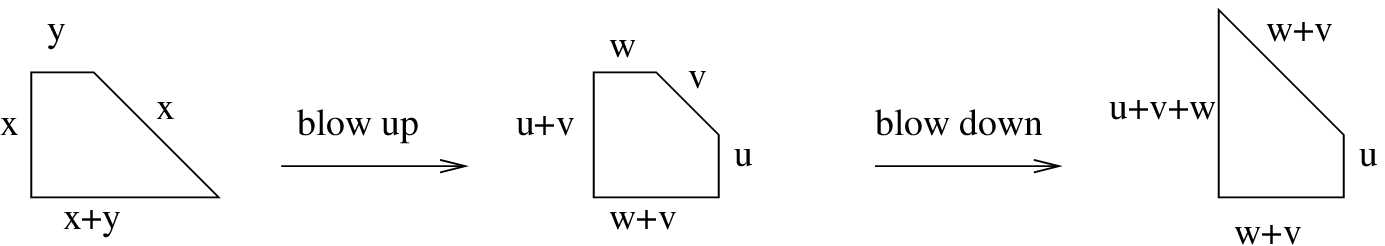}}
       \caption{Reduced Spaces}\label{labelname:4}
       \end{figure}

        Similarly, the Euler classes of the circle bundles change:
          $$-y+u+w=u,$$
         where $\beta^*y=w$.
        
         For $g>0$, since there is no exceptional sphere to blow down in the reduced space $E_{\Sigma_g}$ at a regular value above the minimum, 
         the moment map has to cross
         the index 2 point first. The reduced space changes by blowing up a point in $E_{\Sigma_g}$, then blowing down an exceptional
         sphere. Since after blow down, $M_{red}$ must be $E_{\Sigma_g}$, the  exceptional sphere obtained by blowing up must be
         blown down. Analyze it as in the above (b), we will see that the circle bundles change doesn't satisfy the right formula. So
         $g>0$ can't be true.
        \end{proof} 
        
      About (vii), we can prove the following lemmas:
      
       \begin{lemma} \label{lem9} 
          Assume the fixed point sets are 3 surfaces and there is no twist. Let $g$ and $g_1$ be the genus of the minimum (maximum) and the index 2 surface respectively, and let $b_{min}=n$. Then $b_{max}=-n-2(1+g_1-2g)$;  $\eta=(1+g_1-2g)x+2y$ for even
          $n$, and $\eta=(2+g_1-2g)x+2y$ for odd $n$. ($\eta$ is the dual class of the index 2 surface in $M_{red}$.)
          \end{lemma} 
       \begin{lemma} \label{lem10}
          Assume the fixed point sets are 3 surfaces and there is a twist. Then the surfaces are all spheres,
          and either $b_{min}=0$, $b_{max}=n'=2k'$  and $\eta=x+(1-\frac{n'}{2})y$  or
          $b_{min}=n\ne 0$ ($n$ is even), then $b_{max}=0$ and $\eta=(1-\frac{n}{2})x+y$.
          \end{lemma}

     Proof of Lemma~\ref{lem9}:
     \begin{proof}
        
        If $b_{min}=n=2k$ ( or $ 2k+1$), $b_{max}=n'=2k'$ ( or  $2k'+1$) are even (or odd), then $M_{red}$ is diffeomorphic to the trivial $S^2$ bundle over the Riemann surface
        $\Sigma_g$, $\Sigma_g\times S^2$ ( or the nontrivial $S^2$ bundle over $\Sigma_g$, $E_{\Sigma_g}$). Let $\eta=cx+dy$ be the dual class
        of the embedded image of the index 2 surface in the reduced space, by (\ref{eq2}), 
        Lemma~\ref{lem6} and Lemma~\ref{lem6'}, we have
        \begin{equation} \label{eqcircle}
        kx-y+(cx+dy)=-k'x+y.
        \end{equation} 
       So  $\qquad     d=2, \quad  c=-k'-k.$

          The genus $g_1$ surface is embedded in $M_{red}$, the adjunction formula gives

          $$\int_{\Sigma_g\times S^2}(cx+dy)^2-\int_{\Sigma_g\times S^2}(cx+dy)((2-2g)x+2y)+2=2g_1,$$

          $$\left( \hbox{or} \qquad  \int_{E_{\Sigma_g}}(cx+dy)^2-\int_{E_{\Sigma_g}}(cx+dy)((3-2g)x+2y)+2=2g_1, \right)$$\\
        i.e., 
          \begin{equation} \label{eqadj}
              cd-c-(1-g)d+1=g_1.   
          \end{equation}     
         $$( \hbox{or} \qquad 2cd-d^2-2c-(3-2g)d+2d+2=2g_1.)$$

          Since $d=2$, so $c=1+g_1-2g$  ( or   $c=2+g_1-2g$). Hence
          
           $\qquad \eta=(1+g_1-2g)x+2y$  ( or  $\eta=(2+g_1-2g)x+2y$).

          From $k'=-k-c$, we get $n'=-n-2(1+g_1-2g)$. ( Same for odd case.)\\     
      
        \end{proof}

        Proof of Lemma~\ref{lem10}:
        \begin{proof}
        We have seen that in this case, $g=0$ and $b_{min}=n=2k,  b_{max}=n'=2k'$ are even, and $M_{red}\simeq S^2\times S^2$.
         Note also that there is a twist. So (\ref{eqcircle}) becomes
           $$kx-y+(cx+dy)=-k'y+x,$$
        i.e.,
           $$c=1-k, \quad d=1-k'.$$
            In (\ref{eqadj}), take $g=0$, then
           $$cd-c-d+1=g_1.$$
            
            Combine the above equations, we have $kk'=g_1.$
             If $k=0$ (or $k'=0$), then $g_1=0$ (i.e., if $b_{min}=0$ or $b_{max}=0$, then the fixed point sets are 3 spheres). Then
             $\eta=x+(1-k')y$.
             If $k\ne 0$, then $k'=\frac{g_1}{k}$. So $\eta=(1-\frac{n}{2})x+(1-\frac{2g_1}{n})y$. 
          
          We will compute the cohomology class of the reduced symplectic form at the index 2 critical level and make sure this class
          does contain a symplectic form. We will use the fact that the index 2 surface is a symplectic submanifold in the reduced space,
          and the fact that the size of the maximum is positive. We will see that  $g_1>0$ can not be true. So $g_1$=0.

   In detail, assume $\alpha_0$ is the size of the minimum, $t_0$ is the distance between the minimum and the index 2 surface, and $t_1$ is the distance 
   between the index 2 surface and the maximum. For convenience, let's use the following picture 
   (the left hand side data are the Euler classes of the principal circle bundles over $M_{red}$).

                 \begin{figure}[h!]
    \scalebox{.60}{\includegraphics{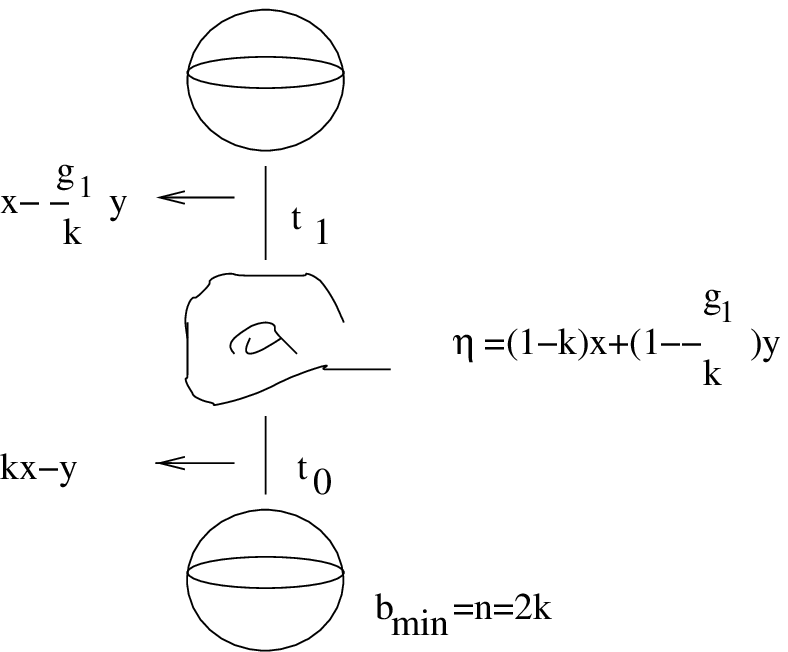}}
    \end{figure}

      By Duistermmat-Heckman formula, when $0<t\leq t_0$, we have $[\omega_t]=(\alpha_0-kt)x+ty$. So we need 
           $$\alpha_0>kt_0$$ 
      to make the classes contain symplectic forms.

      We also need
      $$\int_Z\omega_{t_0}=\alpha_0+(1-2k)t_0+g_1t_0-\frac{g_1}{k}\alpha_0>0.$$

      When $0\leq t\leq t_1$, we have  $[\omega_t]=(\alpha_0-kt_0-t)x+(t_0+\frac{g_1}{k}t)y$.
     At $t=t_1$, we need $\alpha_0-kt_0-t_1=0$, and  $t_0+\frac{g_1}{k}t_1>0.$
     So 
     $$t_1=\alpha_0-kt_0,   \hbox{and} \quad \frac{g_1}{k}\alpha_0-g_1t_0+t_0>0.$$
    In summary, we need
    $$\alpha_0>kt_0,    \hbox{and} \quad \alpha_0+(2-n)t_0>\frac{g_1}{k}\alpha_0-g_1t_0+t_0>0.$$

   Now assume that $g_1>0$, the second part of the second inequality says $k>0$. The second inequality is
      $$\alpha_0+t_0+g_1t_0-nt_0-\frac{g_1}{k}\alpha_0>0.$$
      Since $g_1>0, k>0$ and $\frac{g_1}{k}$ should be an integer, so we can assume $\frac{g_1}{k}=m\geq 1$.
  The last inequality implies 
      $$t_0(1-k)+(1-m)(\alpha_0-kt_0)>0.$$
      According to the choice, this is not possible. Therefore $g_1=0$.
  \end{proof}

      This completes the proof of Theorem 1.
      
     \section{Proof of Proposition 1}

     Let  $S$ be a fixed surface of index 2, and $S\subset \phi^{-1}(c)$ be the only fixed point set on this level. The normal bundle of $S$ in $M$ splits into two
     line bundles: the positive normal bundle $E_+$ and the negative normal bundle $E_-$. By [Mc], the normal bundle of $i(S)=Z$
     in $M_{red}$ is $E_+\otimes E_-$. So if $b_+$ and $b_- $ are the Euler numbers of $E_+$ and $ E_-$ respectively and $\eta$ is the dual class of $Z$ in $M_{red}$,
     then
    \begin{equation} \label{eq+}
         b_+ +b_-=\int_{M_{red}}\eta ^2.  
    \end{equation}

       Let's look at the symplectic cutting ([L]) in a neighborhood of the critical level $c$, we get a 6-dimensional space $M_{c-\epsilon}^{c+\epsilon}$
       whose interior is diffeomorphic to $\phi^{-1}(c-\epsilon, c+\epsilon)$, whose boundary is $M_{c-\epsilon}\cup M_{c+\epsilon}$.
       It has 3 fixed point set components: $M_{c-\epsilon}$ as the minimum with normal bundle $P_{c-\epsilon}$, the index 2 surface $S$ with normal bundles $E_+$ and $E_-$, 
       and $M_{c+\epsilon}$ as the maximum with normal bundle $-P_{c+\epsilon}$.

       On this space, integrate 1:

            $$\int_{M_{c-\epsilon}}\frac{1}{\lambda+e(P_{c-\epsilon})} + \int_{M_{c+\epsilon}}\frac{1}{-\lambda-e(P_{c+\epsilon})} + \int_{S}\frac{1}{(\lambda+b_+u)(-\lambda+b_-u)}=0.$$
      This gives
       
    \begin{equation} \label{eq-}
        b_+- b_-=\int_{M_{c+\epsilon}}e^2(P_{c+\epsilon})-\int_{M_{c-\epsilon}}e^2(P_{c-\epsilon}).
    \end{equation} \\
    \begin{remark}
     If there are other fixed point sets on the same level as $S$, then (\ref{eq-}) is different, and
    $e(P_{c+\epsilon})$ and $e(P_{c-\epsilon})$ will change accordingly.
    \end{remark}

    Use (\ref{eq2}), (\ref{eq6}) and (\ref{eq7}) if necessary, it is easy to know $e(P_{c-\epsilon})$ and  $e(P_{c+\epsilon})$ for each critical level
    $c$ on  which the index 2 surface is.

     Then apply (\ref{eq+}) and (\ref{eq-}), we can obtain:

     For manifold of type (1): $b_+=b_-=2$.

     For manifold of type (2): $ b_-^I=1, b_+^I=0, b_-^{\Pi}=0, b_+^{\Pi}=1$, where $I$ and $\Pi$ denote the first and second index 2 sphere respectively.

     For manifold of type (3): $b_+=1,  b_-=1-n$, where $b_{min}=n=2k+1$.

     For manifold of type (6a): $b_+=n+3(1+g_1-2g), b_-=-n+(1+g_1-2g)$.

     For manifold of type (6b): if $b_{min}=0$, then
     $b_+=1-n', b_-=1$; if $b_{min}=n\neq 0$, then $b_+=1, b_-=1-n$.

   \section{Proof of Theorem 2}

      Let $S^1$ act on a compact symplectic manifold $M$ with moment map $\phi$.  Kirwan ([K])
      shows that the inclusion of the fixed point sets to the manifold $ M^{S^1}\hookrightarrow M$ induces an injection 
      $H_{S^1}^\ast(M,\mathbb{Q})\longrightarrow H_{S^1}^\ast(M^{S^1},\mathbb{Q})$.
      If the cohomology of the fixed point sets has no $p$ torsion for any prime $p$, or the action is semi-free, Tolman and Weitsman  ([TW2]) show
      that this is also true for integral cohomology.\\

      So  we want to determine $H_{S^1}^\ast(M,\mathbb{Z})|_{M^{S^1}}$.\\

    To compute $H_{S^1}^*(M,\mathbb{Z})|_{M^{S^1}}$, let's consider the equivariant Morse stratification. Put an $S^1$ invariant Hermitian metric on $M$ and
    look at the flow of the gradient vector field. The  unstable manifolds corresponding to the fixed point sets form a Morse stratification of $M$.
    The Thom classes of the stratification is a basis of $H^*(M,\mathbb{Z})$ (this is also true for the negative gradient flow). Since the Morse stratification is invariant under the circle action, we 
     have the equivariant Morse stratification, and the corresponding equivariant Thom classes give a basis of $H^*_{S^1}(M,\mathbb{Z})$ as a $H^*_{S^1}(pt)$
     module. For each fixed point set component $F$ with index $i$, there
     is a  class $\alpha\in H^*_{S^1}(M,\mathbb{Z})$ such that $\alpha|_F$ is the equivariant Thom  class of the  negative normal bundle of $F$, and the restrictions 
     of $\alpha$ to the fixed point sets with lower indices are zero. The restriction of this class to the fixed point set with the same index
     depends, but we can make a proper choice. The equivariant integration formula  allows us to obtain this and the restrictions of this class
     to the fixed point sets with higher indices. So we can completely determine $H^*_{S^1}(M,\mathbb{Z})|_{M^{S^1}}$ (See [K] for general idea about the basis).\\

     In detail, we will give this basis explicitly for each of the manifolds we obtained.

     Proof of Theorem 2:
     \begin{proof}
     
      Let $\lambda\in H_{S^1}^2(pt)$ and $u\in H^2(\hbox{surface})$ be generators.\\

     For manifold of type (1), let's label the minimum by $F_1$, the sphere by $F_2$, and the maximum by $F_3$.

       The basis of $H_{S^1}^{\ast}(M,\mathbb{Z})$ given by [K] is:\\
       $\alpha_1 \in H_{S^1}^0(M, \mathbb{Z})$, such that $\alpha_1 \equiv 1.$\\
       $\alpha_2 \in H_{S^1}^2(M, \mathbb{Z})$, such that $\alpha_2|_{F_2}=-\lambda+b_-u=-\lambda+2u$, and  $\alpha_2|_{F_1}=0$. $\alpha_2|_{F_3}$ needs to be determined.\\
       $\alpha'_2 \in H_{S^1}^4(M, \mathbb{Z})$, such that $\alpha'_2|_{F_2}=-\lambda u$, and $\alpha'_2|_{F_1}=0$. $\alpha'_2|_{F_3}$ needs to be determined.\\
       $\alpha_3 \in H_{S^1}^6(M, \mathbb{Z})$, such that $ \alpha_3|_{F_3}=-\lambda^3$, and $\alpha_3|_{F_1}=\alpha_3|_{F_2}=0.$\\
       Now use integration formula to compute the unknowns. ( the denominators are always the equivariant Euler classes of the normal bundles of the fixed
       point sets.)
       
       Integrate $\alpha'_2$ on $M$:
       
          $$ \frac{\alpha'_2|_{F_3}}{-\lambda^3}+\int_{F_2}\frac{-\lambda u}{(\lambda+b_+u)(-\lambda+b_-u)}+0=0.$$
     We get $\alpha'_2|_{F_3}=\lambda^2.$
                     
       Integrate $\alpha_2$ on $M$:
       
          $$\frac{\alpha_2|_{F_3}}{-\lambda^3}+\int_{F_2}\frac{-\lambda+b_-u}{(\lambda+b_+u)(-\lambda+b_-u)}+0=0.$$
    We get           $\alpha_2|_{F_3}=-2\lambda.$  (Note that $b_+=2$.)

           So the equivariant cohomology ring is uniquely determined.\\

      For manifold of type (2), let's label the minimum by $F_1$, the first sphere by $F_2$, the second sphere by $F_3$, and the maximum by $F_4$.

      The basis of $H_{S^1}^{\ast}(M,\mathbb{Z})$  is:\\
      $\alpha_1 \in H_{S^1}^0(M,\mathbb{Z})$, such that $\alpha_1\equiv 1.$\\
      $\alpha_2\in H_{S^1}^2(M, \mathbb{Z})$, such that $\alpha_2|_{F_2}=-\lambda+b_-^Iu=-\lambda+u$, and $\alpha_2|_{F_1}=0.$\\
      $\alpha_3 \in H_{S^1}^2(M, \mathbb{Z})$, such that $\alpha_3|_{F_3}=-\lambda+b_-^{\Pi}u=-\lambda$, and $\alpha_3|_{F_1,F_2}=0.$\\
      $\alpha'_2 \in H_{S^1}^4(M, \mathbb{Z})$, such that $\alpha'_2|_{F_2}=-\lambda u$, and $\alpha'_2|_{F_1}=\alpha'_2|_{F_3}=0.$\\
      $\alpha'_3 \in H_{S^1}^4(M, \mathbb{Z})$, such that $\alpha'_3|_{F_3}=-\lambda u$, and $\alpha'_3|_{F_1}=\alpha'_3|_{F_2}=0.$\\
      $\alpha_4 \in H_{S^1}^6(M, \mathbb{Z})$, such that $\alpha_4|_{F_4}=-\lambda^3$, and $\alpha_4|_{F_1, F_2, F_3}=0.$\\
      
      Integrate $\alpha'_3$ on $M$, we get
      $$\alpha'_3|_{F_4}=\lambda^2.$$
          
      Integrate $\alpha'_2$ on $M$, we get
      $$\alpha'_2|_{F_4}=\lambda^2.$$
           
      Integrate $\alpha_3$ on $M$, we get
      $$\alpha_3|_{F_4}=-\lambda.$$
           
           Assume  $\alpha_2|_{F_3}=cu, \alpha_2|_{F_4}=d\lambda$. 
      
      Integrate $\alpha_2$ on $M$, we get
      $$-d-c=0.$$
           
      Integrate $\alpha_2 \cdot c_1(TM)$, we get
      
         $\frac{d\lambda\cdot (-3\lambda)}{-\lambda^3}+\int_{F_3}\frac{cu\cdot (2u+(b_+^{\Pi}+b_-^{\Pi})u)}{(\lambda+u)(-\lambda+0)}+\int_{F_2}\frac{(-\lambda+u)(2u+(b_+^I+b_-^I)u)}{(\lambda+0)(-\lambda+u)}+0=0.$
         
         We get $d=-1, c=1.$ Hence, $\alpha_2|_{F_4}=-\lambda,  \alpha_2|_{F_3}=u.$

       The equivariant cohomology ring is uniquely determined.\\

       For manifolds of type (3),  let's label the critical sets as: $F_1$ for the minimum, $F_2$ for the index 2 sphere, $F_3$ for the index 4
       point, and $F_4$ for the maximum.

       We have the following basis for $H^*_{S^1}(M, \mathbb{Z})$:\\
       $\alpha_1\in H^0_{S^1}(M, \mathbb{Z})$ such that  $\alpha_1\equiv 1.$\\
       $\alpha'_1\in H^2_{S^1}(M, \mathbb{Z})$ such that $\alpha'_1|_{F_1}=u$. we need to compute the restrictions of $\alpha'_1$ to the other
       fixed point sets.\\
       $\alpha_2\in H^2_{S^1}(M, \mathbb{Z})$ such that  $\alpha_2|_{F_2}=-\lambda+b_-u$, and $\alpha_2|_{F_1}=0$. The restrictions of $\alpha_2$ to the
       other fixed point sets need to be determined.\\
       $\alpha'_2 \in H^4_{S^1}(M, \mathbb{Z})$  such that  $\alpha'_2|_{F_2}=-\lambda u$,  and $\alpha'_2|_{F_1, F_3}=0$.  $\alpha'_2|F_4$ needs to be determined.\\
       $\alpha_3\in H^4_{S^1}(M, \mathbb{Z})$  such that  $\alpha_3|_{F_3}=\lambda^2$, and $\alpha_3|_{F_1, F_2}=0$.  $\alpha_3|_{F_4}$ needs to be determined.\\
       $\alpha_4\in H^6_{S^1}(M, \mathbb{Z})$  such that  $\alpha_4|_{F_4}=-\lambda^3$, and $\alpha_4|_{F_1, F_2, F_3}=0.$

       Integrate $\alpha'_2$:

        $$\frac{\alpha'_2|_{F_4}}{-\lambda^3}+\int_{F_2}\frac{-\lambda u}{(\lambda+b_+u)(-\lambda+b_-u)}=0.$$
       We get
       $$\alpha'_2|_{F_4}=\lambda^2.$$
        
        Integrate $\alpha_3$:

         $$\frac{\alpha_3|_{F_4}}{-\lambda^3}+\frac{\lambda^2}{\lambda^3}=0.$$
      We get
      $$\alpha_3|_{F_4}=\lambda^2.$$

          We can assume  $\alpha'_1|_{F_4}=a\lambda,  \alpha'_1|_{F_3}=b\lambda$, and  $\alpha'_1|_{F_2}=cu$.

          Integrate $\alpha'_1$:

          $$\frac{a\lambda}{-\lambda^3}+\frac{b\lambda}{\lambda^3}+\int_{F_2}\frac{cu}{(\lambda+b_+u)(-\lambda+b_-u)}+\int_{F_1}\frac{u}{\lambda^2+b_{min}\lambda u}=0.$$
         We get
         $$-a+b-c+1=0.$$

          Integrate $\alpha'_1\cdot c_1(TM)$:

          $$\frac{a\lambda\cdot (-3\lambda)}{-\lambda^3}+\frac{b\lambda\cdot (-\lambda)}{\lambda^3}+\int_{F_2}\frac{cu\cdot (2+b_++b_-)u}{(\lambda+b_+u)(-\lambda+b_-u)}$$

            $$+\int_{F_1}\frac{u(2\lambda+2u+b_{min}u)}{\lambda^2+b_{min}\lambda u}=0.$$
        We get
        $$3a-b+2=0.$$

        Integrate ${\alpha'_1}^2$:
        
       $$\frac{a^2\lambda^2}{-\lambda^3}+\frac{b^2\lambda^2}{\lambda^3}+0+0=0.$$
      We get
          $$a^2=b^2.$$

        Combine the above 3 equations, we can get the only integer solution:
               $$a=b=-1,\quad  c=1.$$

        We can assume $\alpha_2|_{F_4}=d\lambda$, and $\alpha_2|_{F_3}=e\lambda$. Similar to the above, integrate $\alpha_2$ and $\alpha_2\cdot c_1(TM)$,
        we obtain 
        $$d=-\frac{2+b_-}{2},  \quad  e=-\frac{b_-}{2}.$$
           
        We can also integrate $\alpha^2_2$ and $\alpha_2\cdot \alpha'_1$ to get another 2 equations, one can check that the above solutions
        will satisfy these 2 equations.

           So we see that the fixed point data determines the equivariant cohomology ring.\\

        Due to similar computations, we will omit some detailed calculations  for manifolds of type (4) and (5).

     For manifold of type (4), let's call the maximum $F_1$, and the minimum $F_2$.
           Note that $b_{max}=b_{min}=2$.
         
          The basis is:\\          
          $\alpha_2\in H^0_{S^1}(M, \mathbb{Z})$ such that $\alpha_2\equiv 1$.\\
          $\alpha'_2\in H^2_{S^1}(M, \mathbb{Z})$ such that $\alpha'_2|_{F_2}=u$, and $\alpha'_2|_{F_1}=a\lambda+bu$, where $a, b$ need
          to be determined.\\
           $\alpha_1\in H^4_{S^1}(M, \mathbb{Z})$ such that $\alpha_1|_{F_1}=\lambda^2-b_{max}\lambda u=\lambda^2-2\lambda u$, and $\alpha_1|_{F_2}=0.$\\
           $\alpha'_1\in H^6_{S^1}(M, \mathbb{Z})$ such that $\alpha'_1|_{F_1}=\lambda^2u$, and $\alpha'_1|_{F_2}=0.$

            Integrate $\alpha'_2$ and $\alpha'_2\cdot c_1(TM)$, we obtain
           $$\alpha'_2|_{F_1}=-\lambda+u.$$
   The equivariant cohomology ring is uniquely determined.\\

        For manifold of type (5), let's label the minimum by $F_1$, the index 2 point by $F_2$, the index 4 point by $F_3$, and the maximum by $F_4$.
          
           Note that $b_{max}=b_{min}=1$.

            The basis is:\\
           $\alpha_1\in H^0_{S^1}(M, \mathbb{Z})$, such that $\alpha_1\equiv 1.$\\
           $\alpha'_1\in H^2_{S^1}(M, \mathbb{Z})$, such that $\alpha'_1|_{F_1}=u, \alpha'_1|_{F_2}=0$, $\alpha'_1|_{F_3}=a\lambda,$
         and  $\alpha'_1|_{F_4}=b\lambda+cu$.\\
           $\alpha_2 \in H^2_{S^1}(M, \mathbb{Z})$, such that $\alpha_2|_{F_2}=-\lambda, \alpha_2|_{F_1}=0, \alpha_2|_{F_3}=d\lambda,$
         and  $\alpha_2|_{F_4}=e\lambda+fu.$\\
           $\alpha_3 \in H^4_{S^1}(M, \mathbb{Z})$, such that $\alpha_3|_{F_3}=\lambda^2,  \alpha_3|_{F_2, F_1}=0$, and  $\alpha_3|_{F_4}=g\lambda u$.\\
           $\alpha_4 \in H^4_{S^1}(M, \mathbb{Z})$, such that $\alpha_4|_{F_4}=\lambda^2-b_{max}\lambda u=\lambda^2-\lambda u$, and $\alpha_4|_{F_1, F_2, F_3}=0.$\\
           $\alpha'_4\in H^6_{S^1}(M, \mathbb{Z})$, such that $\alpha'_4|_{F_4}=\lambda^2u$, and $\alpha'_4|_{F_1, F_2, F_3}=0.$\\
          Where $a, b, c, d, e, f, g$ need to be determined.
          
          Integrate $\alpha_3$, we get
            $$\alpha_3|_{F_4}=-\lambda u.$$
            
          Integrate $\alpha'_1$, $\alpha'_1\cdot c_1(TM)$ and $(\alpha'_1)^2$, we get
            $$\alpha'_1|_{F_4}=-\lambda+u, \quad \alpha'_1|_{F_3}=-\lambda.$$ 
            
            Similarly, integrate $\alpha_2$, $\alpha_2\cdot c_1(TM)$ and $(\alpha_2)^2$, we get
             $$\alpha_2|_{F_4}=-\lambda, \quad \alpha_2|_{F_3}=0.$$
             
             (If we integrate $\alpha_2\cdot \alpha'_1$, we will see that the above solutions satisfy this equation.)
          
           Hence the fixed point data determines the  equivariant cohomology ring.\\

           For manifolds of type (6a) and (6b), let's label the minimum surface by $F_1$, the index 2 surface by $F_2$, and the maximum surface by $F_3$.

          The basis is:\\
           $\alpha_1\in H^0_{S^1}(M, \mathbb{Z})$, such that $\alpha_1\equiv 1.$\\
           $\alpha'_1\in H^2_{S^1}(M, \mathbb{Z})$, such that $\alpha'_1|_{F_1}=u,  \alpha'_1|_{F_2}=\tilde du$, and $\alpha'_1|_{F_3}=e\lambda+fu$, where
           $\tilde d, e, f$ will be determined.\\
           $\alpha_2\in H^2_{S^1}(M, \mathbb{Z})$, such that $\alpha_2|_{F_2}=-\lambda+b_-u, \alpha_2|_{F_1}=0$, and $\alpha_2|_{F_3}=a\lambda+bu$, where
           $a, b$ need to be determined.\\
           $\alpha'_2\in H^4_{S^1}(M, \mathbb{Z})$, such that $\alpha'_2|_{F_2}=-\lambda u, \alpha'_2|_{F_1}=0$, and $\alpha'_2|_{F_3}=c\lambda u$, where
           $c$ needs to be determined.\\ 
           $\alpha_3\in H^4_{S^1}(M, \mathbb{Z})$, such that $\alpha_3|_{F_3}=\lambda^2-b_{max}\lambda u$, and $\alpha_3|_{F_1,F_2}=0.$\\
           $\alpha'_3\in H^6_{S^1}(M, \mathbb{Z})$, such that $\alpha'_3|_{F_3}=\lambda^2u,  \alpha'_3|_{F_1}=0$, and $\alpha'_3|_{F_2}=0.$

           Integrate $\alpha'_2$, we obtain
            $$\alpha'_2|_{F_3}=-\lambda u.$$
          
           Integrate $\alpha'_1$:
           
           $$\int_{F_1}\frac{u}{\lambda^2+b_{min}\lambda u}+\int_{F_2}\frac{\tilde du}{(\lambda+b_+u)(-\lambda+b_-u)}+\int_{F_3}\frac{e\lambda+fu}{\lambda^2-b_{max}\lambda u}=0.$$
         We get

           \begin{equation}\label{eq3.2.1}
             b_{max}e+f-\tilde d+1=0.
           \end{equation}

           Integrate $\alpha'_1\cdot c_1(TM)$:
           
           $$\int_{F_1}\frac{u(2\lambda+(2-2g)u+b_{min}u)}{\lambda^2+b_{min}\lambda u}+\int_{F_2}\frac{\tilde du((2-2g_1)u+(b_++b_-)u)}{(\lambda+b_+u)(-\lambda+b_-u)}$$
           
           $$+\int_{F_3}\frac{(e\lambda+fu)(-2\lambda+(2-2g)u+b_{max}u)}{\lambda^2-b_{max}\lambda u}=0.$$
         We get
           
           \begin{equation}\label{eq3.2.2}
            (2-2g-b_{max})e-2f+2=0.
           \end{equation}

           Integrate $(\alpha'_1)^2$:
           
           $$\int_{F_3}\frac{(e\lambda+fu)^2}{\lambda^2-b_{max}\lambda u}+0+0=0.$$
         We get
           
           \begin{equation}\label{eq3.2.3}
            b_{max}e^2+2ef=0.
           \end{equation}

            Integrate $\alpha'_1\cdot \alpha_2$:
            
            $$\int_{F_3}\frac{(a\lambda+bu)(e\lambda+fu)}{\lambda^2-b_{max}\lambda u}+\int_{F_2}\frac{\tilde du(-\lambda+b_-u)}{(\lambda+b_+u)(-\lambda+b_-u)}+0=0.$$
        We get
            
            \begin{equation}\label{eq3.2.4}
            b_{max}ae+(af+be)+\tilde d=0.
            \end{equation}

            Similarly,     
            integrate $\alpha_2$, we get 
             \begin{equation}\label{eq3.2.5}
             b_{max}a+b-b_+=0.
             \end{equation}

            Integrate $\alpha_2\cdot c_1(TM)$, we get
             \begin{equation}\label{eq3.2.6}
             (2-2g-b_{max})a-2b+2-2g_1+(b_++b_-)=0.
             \end{equation}

            Integrate $\alpha_2^2$, we get 
            \begin{equation}\label{eq3.2.7}
             b_{max}a^2+2ab+(b_++b_-)=0.
             \end{equation}\\
             
      Note that in the above computation, for manifolds of type (6b), $b_{max}=0, g=g_1=0$.

      If we solve the equations directly, we will get multiple solutions. Hence we need the following lemma to finish the proof. 
             \begin{lemma}\label{lem12}
             Assume the fixed point set consists of  3 surfaces and $b_{min}$ is even. Let $F_2$ be the index 2 surface. Let $\alpha'_1\in H^2_{S^1}(M,\mathbb{Z})$
             be the class such that $\alpha'_1|_{\hbox{min}}=u$, $\alpha'_1|_{F_2}=\tilde du$  and 
             $\alpha'_1|_{\hbox{max}}=e\lambda+fu$. 
             Then we always have $\tilde d=d$, where $d$ is the $y$ coefficient of the dual class
             of the index 2 surface  in $M_{red}$. Moreover, for the non-twisted case $e=0$, and for the twisted case  $e=-1$.
             \end{lemma}
             \begin{proof}
             Let's assume $\phi(M)=[0, h]$, and let $\phi(F_2)=r$.

               Consider the following commutative diagram
               \[
                 \begin{array}{ccc}
                 H^*_{S^1}(M)&\longrightarrow  & H^*_{S^1}(\phi^{-1}(\epsilon))\\
                 \downarrow  &                 & \downarrow\\
                 H^*_{S^1}(F_2) &           &          \\
                 \downarrow &                    &                     \\
                 H^*(F_2)&\longleftarrow &      H^*(M_{red})\\ 
                 \end{array}
               \]
               where $\epsilon<r$.
               The two  downward maps on the left clearly map $\alpha'_1$ to $\tilde du$ which is an ordinary cohomology class. Let's denote the composition of the downward map on the right
               with the rightward map by $\kappa$.

              Now consider $\alpha'_1$ in the neighborhood $\phi\leq\epsilon$, according to the choice of $\alpha'_1$, it is an ordinary cohomology
              class in this neighborhood, and $\int_{\hbox{min}}\alpha'_1=1$, so its image in the
              reduced space is the dual class of the fiber $S^2$, which is $x$.
                          $$\kappa(\alpha'_1)=x.$$
                
             Regarding the leftward map, it is induced by the following embedding              
               $$F_2\hookrightarrow M_{red}.$$ 
               So
               $$\int_{F_2}(\tilde du)=\tilde d=\int_{M_{red}}\kappa(\alpha'_1)\eta=\int_{M_{red}}x(cx+dy)=d.$$

             Since there is a natural diffeomorphism between the reduced spaces at different regular values, the image of $\alpha'_1$
             in each reduced space remains the same.

               To prove the second statement, we need to consider another cohomology class $\beta\in H^2_{S^1}(M)$ for the downward flow (or for
               the Morse function $-\phi$), such that $\beta|_{\hbox{max}}=u$. We can similarly prove that for the non-twisted case, $\kappa(\beta)=x$,
               for the twisted case, $\kappa(\beta)=y$.
               
               Integrate $(\alpha'_1\cdot \beta)$ over $\bar{M}_{h-\epsilon}$, (the cut space including the maximum.) (see [L] or [Ka]), we get:\\
               For the non-twisted case,                
                $$\int_{M_{red}}\frac{x\cdot x}{\lambda+e(P_{h-\epsilon})}=0=-\int_{\hbox{max}}\frac{(e\lambda+fu)u}{\lambda^2-n'\lambda u}=-e.$$
               For the twisted case,                 
                $$\int_{M_{red}}\frac{x\cdot y}{\lambda+e(P_{h-\epsilon})}=1=-\int_{\hbox{max}}\frac{(e\lambda+fu)u}{\lambda^2-n'\lambda u}=-e.$$
                
              \end{proof}

               The solution for type (6a):
            $$a=-2,\quad b=-n-(1+g_1-2g), \quad e=0, \quad f=1,\quad \tilde d=2.$$
               
             The solution for type (6b):\\
          For $b_{min}=n=0$, we have
             $$a=-1,\quad b=1,\quad \tilde d=1-\frac{n'}{2},\quad e=-1, \quad f=\frac{n'}{2}.$$
          For $b_{min}=n\neq 0$, we have
             $$a=\frac{n-2}{2},\quad b=1,\quad e=-1, \quad f=0,\quad \tilde d=1.$$
                 \end{proof}

      Therefore, the equivariant cohomology rings of manifolds of type (6a) and type (6b)
    are determined by the fixed point data.

   \section{Proof of Proposition 2}

    \begin{lemma}\label{lem11}
  (See [Li])   Manifolds of type (1), (2), (3), (4), (5) and (6b) are simply connected.  Manifolds of type (6a) with $g=0$ are simply connected, they are not otherwise.
    \end{lemma}

    \begin{lemma}\label{lem12}
     The manifolds in Theorem 1 have no torsion homology.
    \end{lemma}
    \begin{proof}
     This is because the fixed point sets have no torsion homology. See [F].
    \end{proof} 
    \begin{lemma}\label{lem13}
     Let $c_1(TM)$ be the equivariant 1st Chern class of $M$. Let $\alpha_1\equiv 1\in H^*_{S^1}(M,\mathbb{Z})$, and let $\alpha_2$, $\alpha_3$,...
     be a basis of $H^2_{S^1}(M,\mathbb{Z})$. If $c_1(TM)=n_1\lambda\alpha_1+n_2\alpha_2+n_3\alpha_3+...$, and $n_2,n_3,...$ are even integers,
     then $w_2(M)=0$.
    \end{lemma}
   \begin{proof}
  Let $c$ be the ordinary 1st Chern class of the tangent bundle $TM$, then by [MS],  $w_2(M)=c$(mod 2). The lemma follows from this fact.
   \end{proof}

       Assume $M$ is of type (1). Then $c_1(TM)|_{F_1}=3\lambda,  c_1(TM)|_{F_2}=2u+(b_++b_-)u=6u,  c_1(TM)|_{F_3}=-3\lambda$. 
          So  $c_1(TM)=3\lambda\alpha_1+3\alpha_2$. Hence we have $w_2(M)\ne 0$.

     Assume $M$ is of type (2). Then $c_1(TM)=3\lambda\alpha_1+3\alpha_2+3\alpha_3$. So we see $w_2(M)\ne 0$.

     Assume $M$ is of type (3). Then $c_1(TM)|_{F_1}=2\lambda+2u+nu, c_1(TM)|_{F_2}=2u+(b_++b_-)u=2u+(2-n)u, c_1(TM)|_{F_3}=-\lambda-\lambda+\lambda=-\lambda, c_1(TM)|_{F_4}=-3\lambda$.
      So  $c_1(TM)=2\lambda\alpha_1+(2+n)\alpha'_1+2\alpha_2$. Since $n$ is odd, we have $w_2(M)\ne 0$.

     Assume $M$ is of type (4). Then $c_1(TM)=2\lambda\alpha_2+4\alpha'_2$, so $w_2=0$.

     Assume $M$ is of type (5). Then $c_1(TM)=2\lambda\alpha_1+3\alpha'_1+\alpha_2$. So  $w_2\ne 0$.

     Similarly, for $M$ of type (6a), $w_2(M)=0$ when $b_{min}$ is even; for $M$ of type (6b), $w_2=0$.

     Corollary 2 and the lemmas of this section help to check the conditions of Wall's theorem. Proposition 2 follows from there. Since $\mathbb{C}P^3$ with the circle action in the introduction has the same fixed point data as manifold of type (4), so all manifolds of type (4) are diffeomorphic to $\mathbb{C}P^3$. See next section for the toric variety.

    \section{Toric varieties and Proof of Proposition 3}

        Let's recall some basic properties about 6-dimensional toric varieties.

         If $(M^6, \omega)$ has a Hamiltonian $T^3$ action, by Atiyah, Guillemin and Sternberg, the image of the moment map is a 3-dimensional convex polytope 
        in $\mathbb{R}^3$. Conversely, by a theorem of Delzant ([D]),  if we have a convex polytope in $\mathbb{R}^3$, such that for each vertex, the 3 normal
        directions of the 3 faces through the vertex span a $\mathbb{Z}$ basis, then this polytope uniquely determines a 6-dimensional symplectic
        manifold with a $T^3$ Hamiltonian action whose moment map image is exactly the polytope. We call such a polytope a toric variety.\\

         Now if we have a Delzant polytope in $\mathbb{R}^3$, and assume that the 3rd circle of $T=(S^1)^3$ is the circle we are considering acting on 
         $(M^6,\omega)$, then the image of the moment map of the circle action is the image of the projection of the polytope to the 3rd circle
         direction, say, the $z$ direction in $\mathbb{R}^3$ for the coordinate system $xyz$.
         
         Moreover, if we choose the  Delzant polytope in $\mathbb{R}^3$, such that with respect to the 3rd circle action, the action
         is semi-free with the right fixed point sets, then we have a toric variety which corresponds to the manifold.\\

          For this purpose, we need the polytope to satisfy the following conditions:

           About the faces: (a) no horizontal faces (b) if a face has normal direction $(m,n, p)$, then $(m,n)=1.$

           About the edges: if an edge has direction $(m',n',p')$, then either $p'=0$, which corresponds to a fixed sphere for the $S^1$ action,
            or $p'=1$ on which the action is free.

            About the vertices: a vertex which is not on a horizontal edge is an isolated fixed point for the $S^1$ action.\\

    The toric variety corresponding to manifold of type (4) is: \\
       The normal directions of the faces:\\
         face (1) abd, (1,0,0),\\
         face(2)abc, (m,1,0),\\
             face(3)acd, (m,1,1),\\
             face(4)bcd, (-1,0,-1),\\
             where $m$ is an integer (different $m$ gives the same toric variety).
         The picture for $m=2$  is 
     \begin{figure}[h!]
    \rotatebox{270}{\scalebox{.20}{\includegraphics{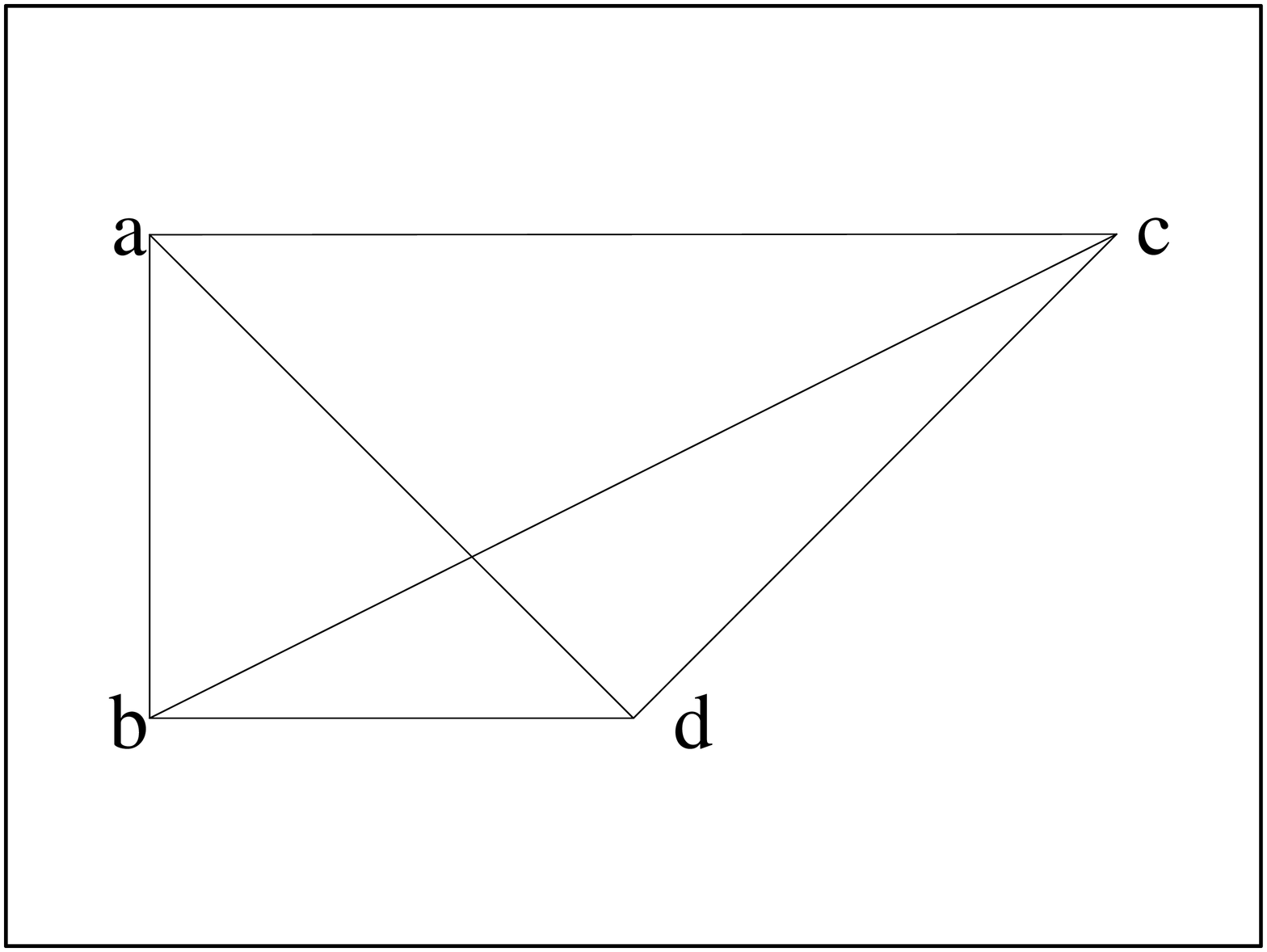}}}
    \end{figure}

     In this picture, the edges bd and ac represent the minimum and the maximum. We saw that the reduced spaces of this manifold are all diffeomorphic 
     to $S^2\times S^2$ (see Lemma~\ref{lem6}). The intersections of this
     picture with horizontal planes are the reduced spaces. It is easy to see that  different homology 
     classes in $S^2\times S^2$ shrink to $0$ when we approach the minimum and the maximum.\\

            The toric variety corresponding to  manifold of type (6b) with $b_{min}=2$ is:\\
            The normal directions of the faces:\\
            face(1)abd: (1,0,0),\\
            face(2)aceb: (m,1,0),\\
            face(3)acgd: (m,1,1),\\
            face(4)cge: (-1,0,0),\\
            face(5)bdge:(-1,0,-1),\\ 
            where $m$ is an integer.
            
            The picture for $m=2$ is 
               \begin{figure}[h!]
    \rotatebox{270}{\scalebox{.20}{\includegraphics{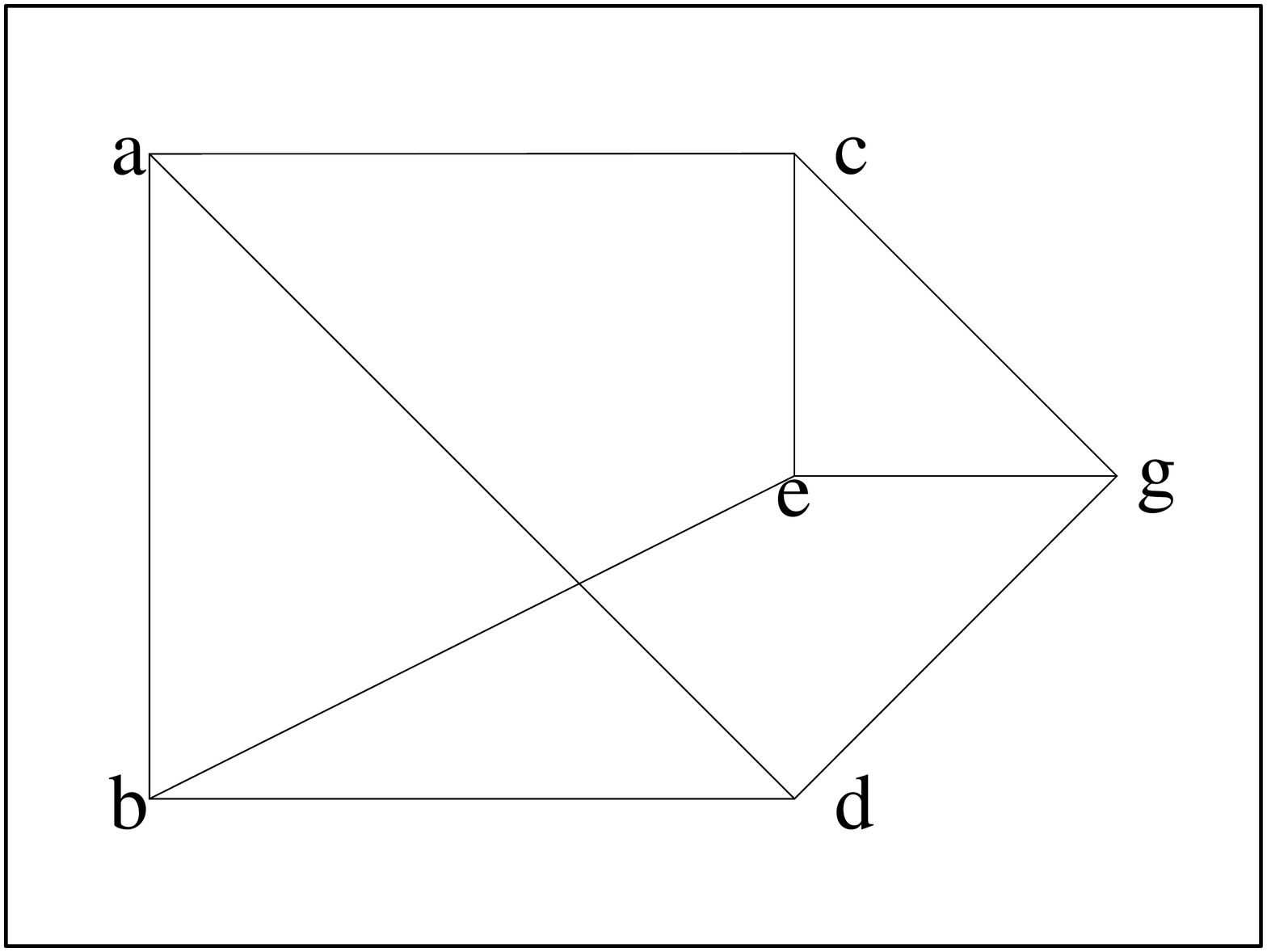}}}
    \end{figure}

    Notice that since  manifold (4) is diffeomorphic to $\mathbb{C}P^3$,  this manifold is  diffeomorphic to $\mathbb{C}P^3$ blown up at a point.\\

       For manifold(s) of type (3),  when $b_{min}=1$ or $b_{min}=3$,
       we have the following toric varieties.

        For $b_{min}=1$ ($b_+=1, b_-=0$), the toric variety  corresponds to $\mathbb{C}P^2\times \mathbb{C}P^1$.

            Normal directions for the faces:\\
                 face (1) abd: (1,0,0),\\
                 face (2) aceb: (m, 1, 0),\\
                 face (3)acgd: (m+1,1,1),\\
                 face (4)cge: (-1, 0, 0),\\
                 face (5)bdge: (-1,0,-1),\\
            where $m$ is an integer. A picture of it when $m=2$ is:
            
    \begin{figure}[h!]
    \rotatebox{270}{{\scalebox{.2}{\includegraphics{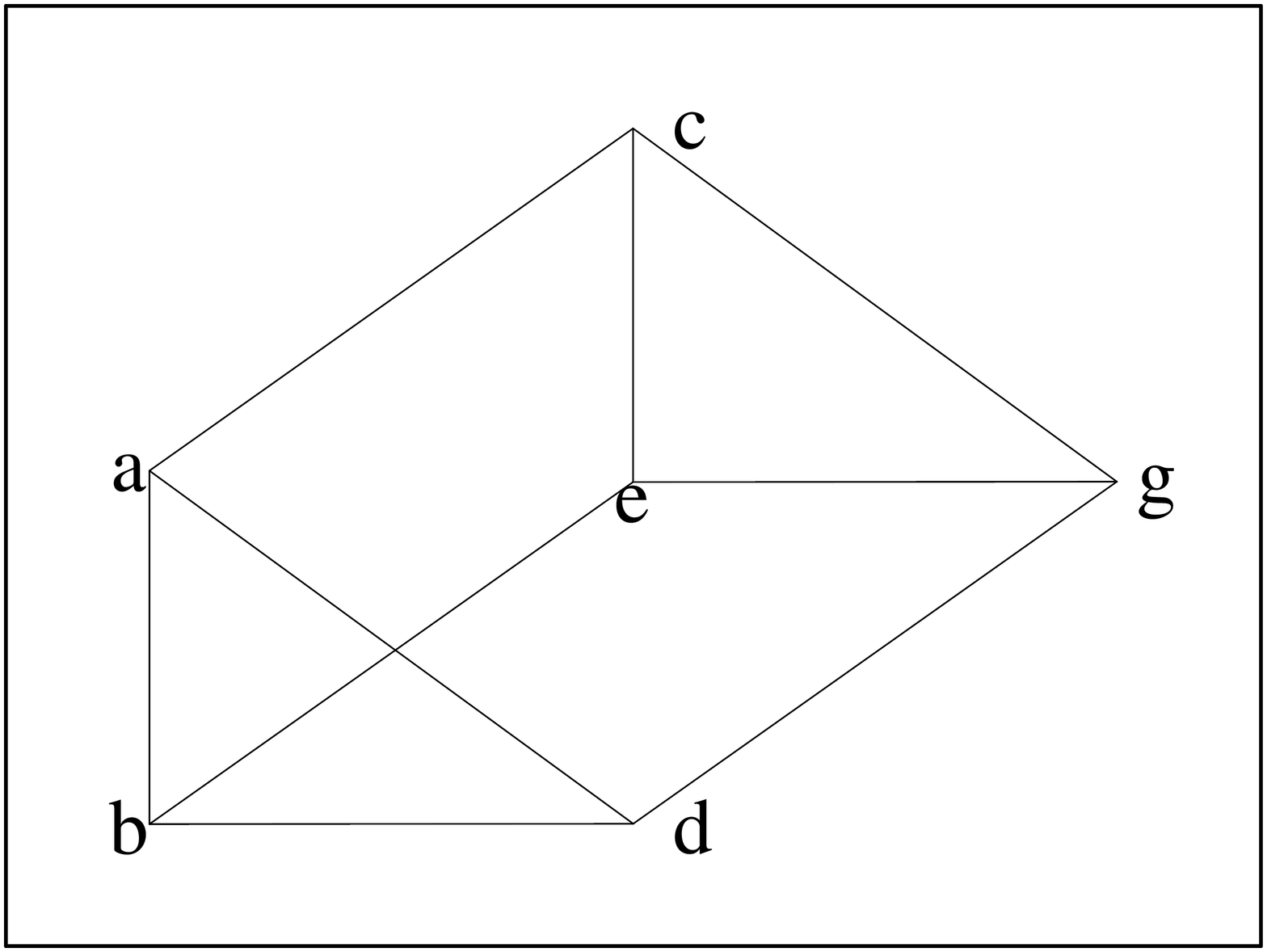}}}}
    \end{figure}

     For $b_{min}=3$ ($b_+=1, b_-=-2$), the toric variety is:\\
     The normal directions of the faces:\\
                face (1)abd: (1,0,0),\\
                face (2)aceb: (m,1,0),\\
                face (3)acgd: (m-1, 1,1),\\
                face (4) cge: (-1, 0,0),\\
                face (5)bdge: (-1, 0,-1),\\
               where $m$ is an integer.

                The picture for $m=2$  is:
          \begin{figure}[h!]
    \rotatebox{270}{{\scalebox{.2}{\includegraphics{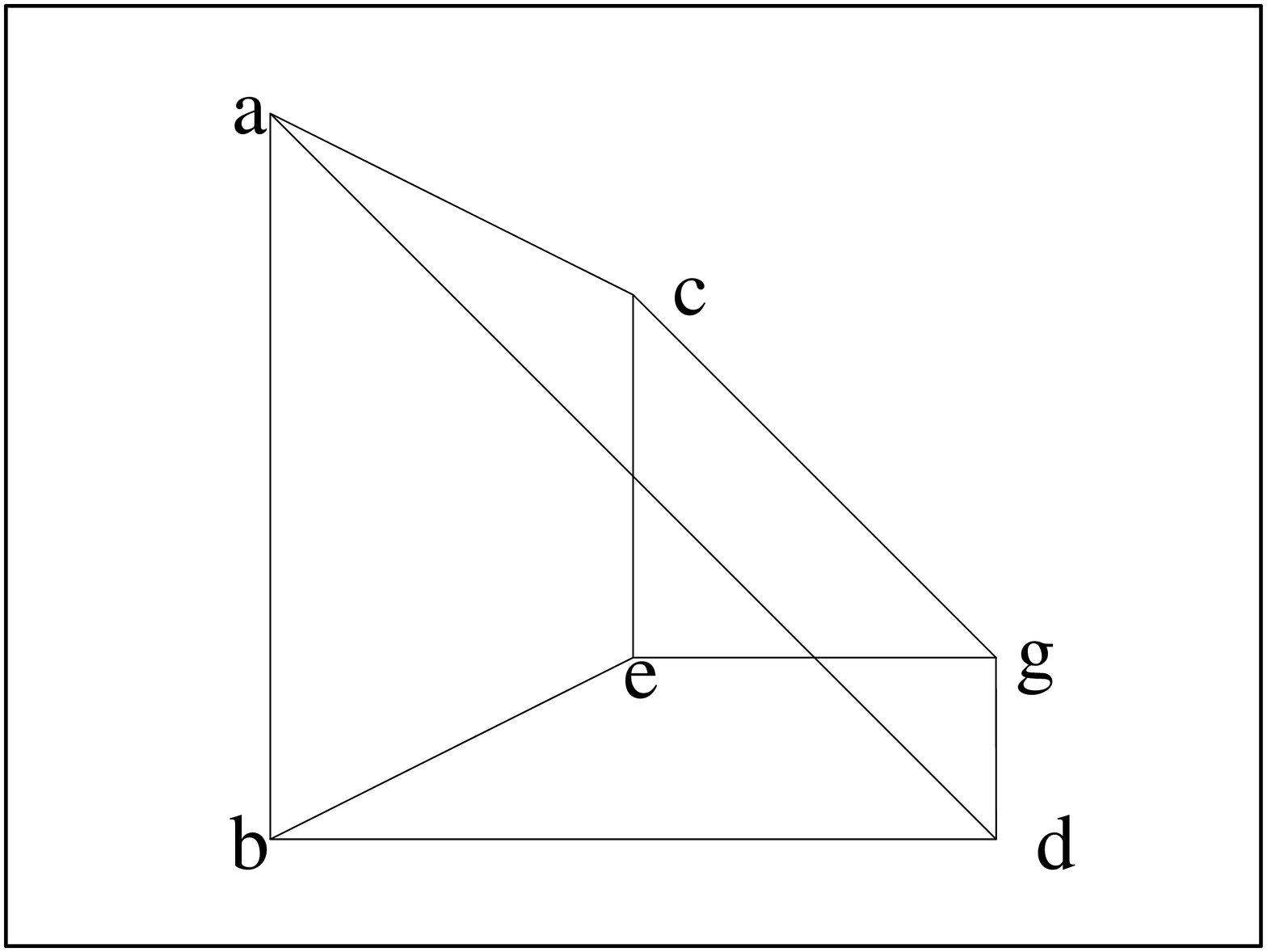}}}}
    \end{figure} 
  \newpage  
    
  Next  we give the two examples  mentioned in Remark~\ref{rem0} in the introduction. Here are two toric varieties which have the same fixed point data but one has no twist and the other one has a twist.

   The example with no twist is the following.

  The normal directions of the faces:\\
  face (1)abcd (1,0,0),\\
  face (2)efgh (-1,0,0),\\
  face (3)bcgf (-1,0,-1),\\
  face (4)abfe (-m,-1,0),\\
  face (5)cghd (m,1,0),\\
  face (6)aehd (1,0,1),\\
  where $m$ is an integer. The picture for $m=2$ is:
      \begin{figure}[h!]
    \rotatebox{360}{{\scalebox{.3}{\includegraphics{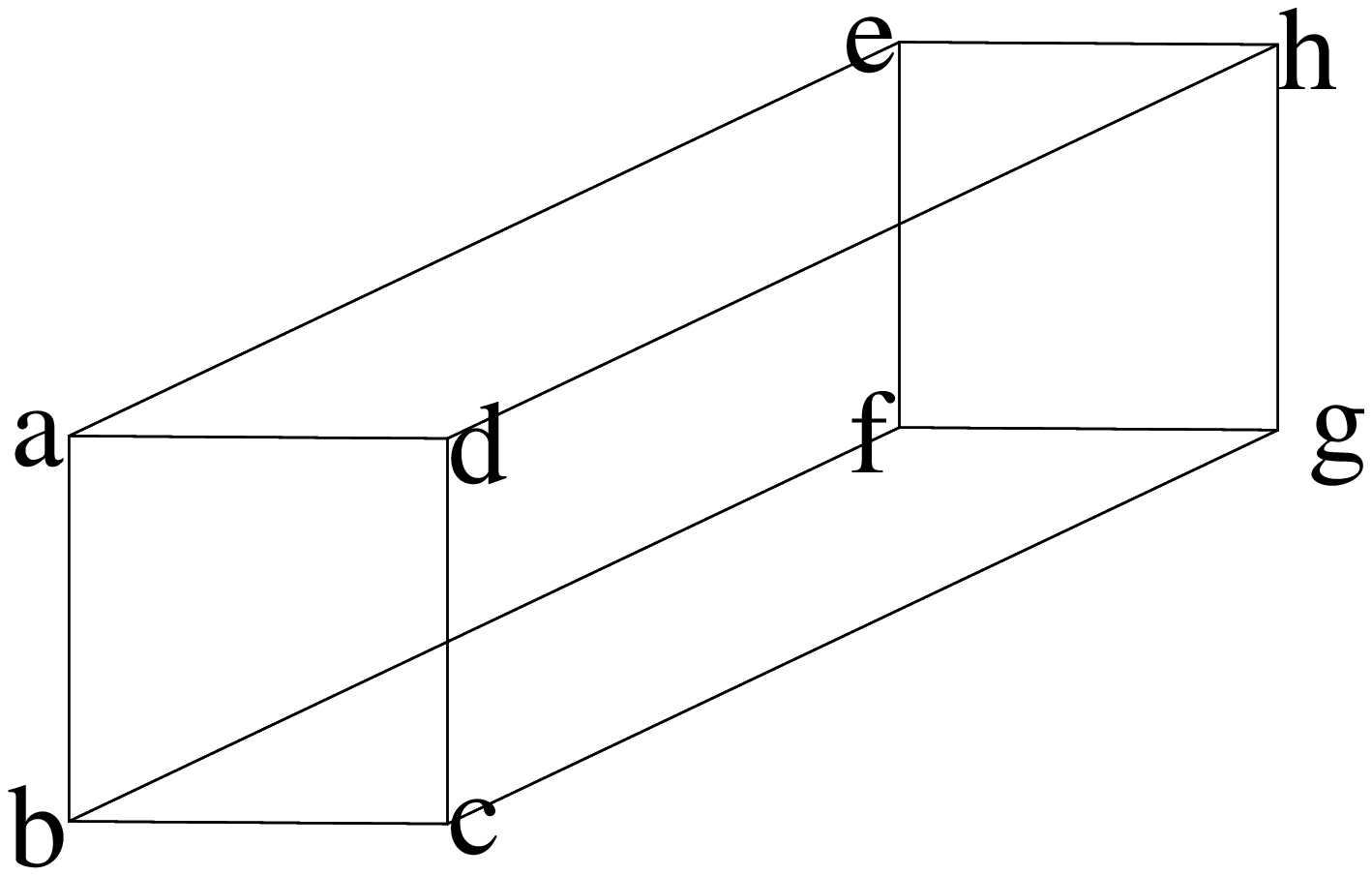}}}}
    \end{figure}

   The example with a twist is the following.

  The normal directions of the faces:\\
  face (1)abcd (1,0,0),\\
  face (2)efgh (-1,0,0),\\
  face (3)bcgf (-1,0,-1),\\
  face (4)abfe (-m,-1,0),\\
  face (5)cghd (m,1,0),\\
  face (6)aehd (m,1,1),\\
  where $m$ is an integer. The picture for $m=2$ is:

     \begin{figure}[h!]
    \rotatebox{360}{{\scalebox{.3}{\includegraphics{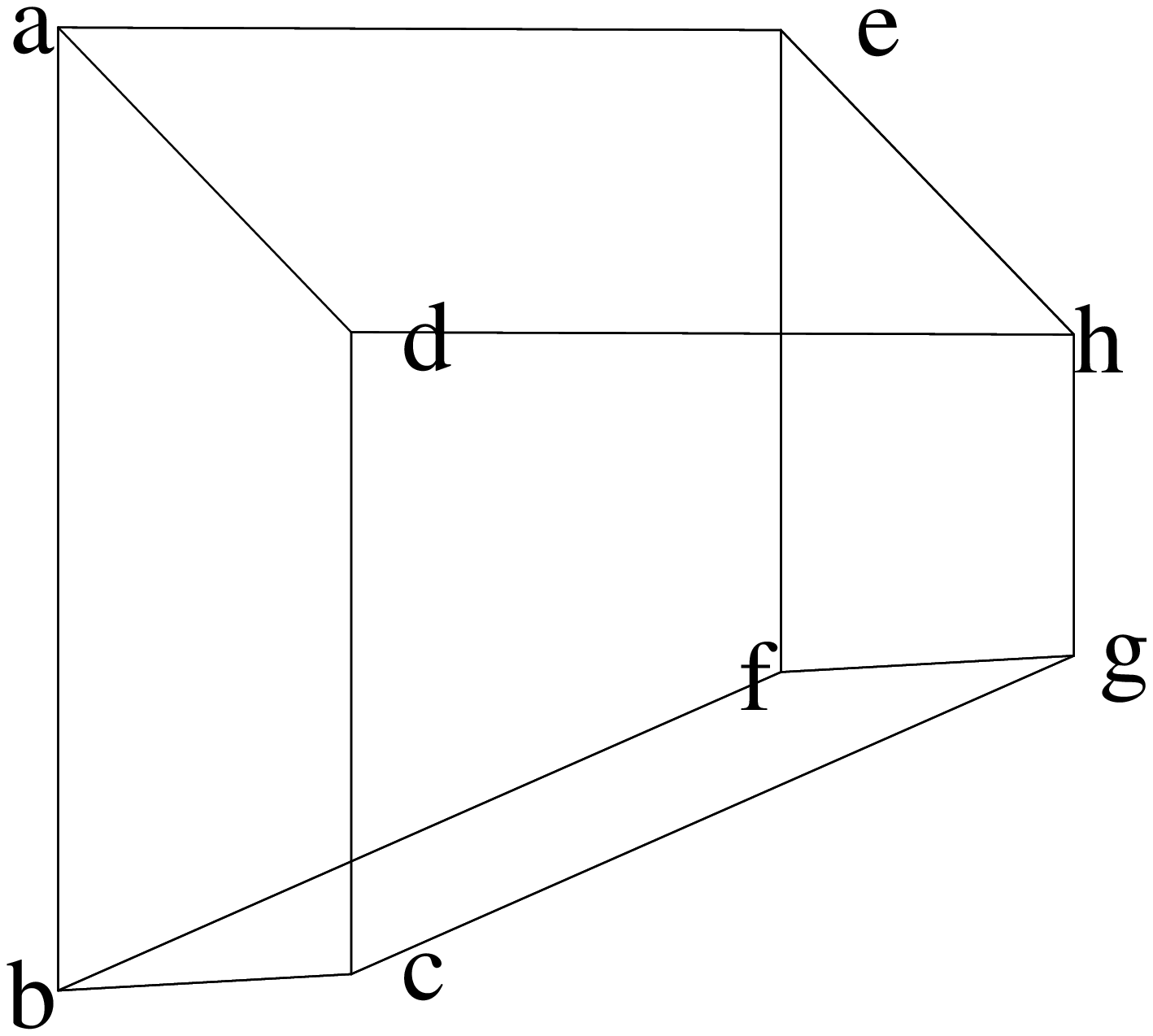}}}}
    \end{figure}

 Notice that these two toric varieties both have 4 fixed spheres (with respect to the circle action) and the Euler numbers of the normal bundles all vanish.

 \end{document}